\documentclass[10pt]{amsart}
\usepackage{amsmath, amssymb, amsthm, graphicx}
\usepackage{epstopdf}

\theoremstyle{plain}
\newtheorem{prop}{Proposition}[section]
\newtheorem{theo}[prop]{Theorem}
\newtheorem{coro}[prop]{Corollary}
\newtheorem{lemm}[prop]{Lemma}

\theoremstyle{definition}

\newtheorem{defi}[prop]{Definition}

\newtheorem{exam}[prop]{Example}
\newtheorem*{rema}{Remark}

\numberwithin{equation}{section}

\renewcommand\aa{a}
\newcommand\act{\cdot}
\newcommand\bb{b}
\newcommand\BB[1]{B_{#1}^{\scriptscriptstyle+}}
\newcommand\Bi{B_\infty}
\newcommand\br[1]{\widehat{#1}}
\newcommand{\cl}[1]{\overline{\vrule height5pt width0pt#1}}
\newcommand\DD[1]{\Delta_{#1}}
\newcommand\divel{\preccurlyeq_{\scriptscriptstyle\! L}}
\newcommand\dual[1]{{}^*\!#1}
\newcommand{\e}{\varepsilon}
\newcommand{\ee}{e}

\newcommand{\etc}{{\it etc.}}
\newcommand\FDyn{F}
\newcommand\ff{f}
\newcommand\flip{\phi}
\newcommand\gcdl{\mathrm{gcd}_{\scriptscriptstyle\! L}}
\newcommand\GDyn{\overline{F}}
\let\ge=\geqslant
\renewcommand\gg{g}
\newcommand\gives{\rightsquigarrow}
\newcommand\HH{\alpha}
\newcommand\id{\mathrm{id}}
\def\ie{{\it i.e.}}
\newcommand\ii{i}
\newcommand{\Int}{\mathbb{Z}}
\newcommand{\inv}{^{-1}}
\newcommand{\jj}{j}
\newcommand{\kk}{k}
\newcommand\lcml{\mathrm{lcm}_{\scriptscriptstyle\! L}}
\let\le=\leqslant
\newcommand{\mm}{m}
\newcommand\Neg[1]{#1^-}
\newcommand{\nn}{n}
\newcommand{\NN}{N}
\newcommand{\perm}{\pi}
\newcommand\Pos[1]{#1^+}
\newcommand\pp{p}
\newcommand{\qq}{q}
\newcommand{\redrr}{\curvearrowright}
\newcommand{\redrl}{\mathrel{\raisebox{5pt}{\rotatebox{180}{$\curvearrowleft$}}}}
\newcommand{\Ree}{\mathbb{R}}
\def\resp{{\it resp.}~}

\newcommand\s{\sigma}
\renewcommand\ss[1]{\sigma_{#1}^{\vrule height5pt width0pt}}
\newcommand\sss[1]{\sigma_{#1}^{-1}}
\newcommand\ssss[1]{\sigma_{#1}^{\pm1}}
\newcommand\tta{\mathtt{a}}
\newcommand\ttA{\mathtt{A}}
\newcommand\ttb{\mathtt{b}}
\newcommand\ttB{\mathtt{B}}
\newcommand\ttc{\mathtt{c}}

\newcommand{\uu}{\ww'}
\newcommand{\vs}{{\it vs.}~}
\newcommand{\vv}{\ww''}
\newcommand{\ww}{w}
\newcommand\xs{s}
\newcommand\xt{t}
\newcommand\xu{u}
\newcommand\xv{v}
\newcommand\xx{x}
\newcommand{\yy}{y}
\newcommand{\zz}{z}

\begin{document}

\hfill\raisebox{20pt}{$\scriptscriptstyle 2007-03$}

\author{Patrick DEHORNOY}
\address{Laboratoire de Math\'ematiques Nicolas Oresme,
UMR 6139 CNRS,  Universit\'e de Caen BP 5186, 14032 Caen,
France}
\email{dehornoy@math.unicaen.fr}
\urladdr{//www.math.unicaen.fr/\!\hbox{$\sim$}dehornoy}

\title{Efficient solutions to the braid isotopy problem}


\thanks{}


\begin{abstract}
We describe the most efficient solutions to the word problem
of Artin's braid group known so  far, \ie, in other
words, the most efficient solutions to the braid isotopy
problem, including the Dynnikov method, which could be especially
suitable for cryptographical applications. Most results appear in
literature; however, some results about the greedy normal form 
and the symmetric normal form and their connection with grid 
diagrams may have never been stated explicitly.  
\end{abstract}

\maketitle

Because they are both not too simple and not too complicated,
Artin's braid groups~$B_n$ have been and remain one of the most
natural and promising platform groups for non-commutative
group-based cryptography~\cite{AAG, KLC, Dgw}. More precisely,
braid groups are not too simple in that they lead to problems with
presumably difficult instances, typically the conjugacy problem and
the related conjugacy  and multiple conjugacy search problems, and
they are not too complicated in that there exist efficient solutions to
the word problem, a preliminary requirement when one aims at
practically computing in a group. It turns out that, since the
founding paper~\cite{Art} appeared in~1947, the word problem of
braid groups---which is also the braid isotopy problem---has
received a number of solutions: braid groups might even be the
groups for which the number of known solutions to the word
problem is currently the highest one. 

In this paper, we review some of these solutions, namely those that, at
the moment, appear as the most efficient ones for practical
implementation, and, therefore, the most promising ones for
cryptographical applications. What makes the subject specially
interesting is that these solutions relie on deep underlying
structures that explain their efficiency. Five solutions are
described, and they come in two families, namely those
based on a normal form, and those that use no normal form. In the
first family, we consider the so-called greedy normal form, both in its
non-symmetric and  symmetric versions. In the second family, we
consider the so-called word redressing method, which, like the greedy
normal forms, has a quadratic complexity, the handle
reduction method, whose complexity remains unknown
but which is very efficient in practice, and Dynnikov's coordinization
method, which relies on an entirely different, geometric approach,
and might turn to be very efficient. In this description, we
deliberately discard lots of alternative solutions which are
intrinsically exponential in complexity.

The paper is organized as follows. In Section~\ref{S:Normal}, after
setting the background, we describe the greedy normal form and the
symmetric normal form, without providing explicit rules to compute
them. All results in this section are standard. In Section~\ref{S:Grid},
we introduce grid diagrams and explain---and prove---how to use
such diagrams to compute the normal forms of
Section~\ref{S:Normal}. Though more or less equivalent to that
of~\cite[Chapter 9]{Eps}, this approach is less standard, and, to the
best of our knowledge, the results may have never been stated in the
form given here.  Finally, in Section~\ref{S:Direct}, we describe the
word redressing, handle reduction, and Dynnikov coordinates
methods. The results here already appeared in literature, but 
Dynnikov's approach, which appears in~\cite[Chapter~7]{Dgr},
has not yet become classical. Also, the formulae of Section~\ref{S:Dynnikov} 
have been optimized to make implementation easy.

\section{Solutions based on a normal form}
\label{S:Normal}

This section deals about solutions to the braid word problem that
consist in defining for each braid~$\xx$ a unique distinguished
representative called the normal form of~$\xx$. When this is done,
in order to compute with braids, it is in practice sufficient to work
with normal representatives. There exist excellent normal forms
for braids, namely those connected with the so-called {\it greedy 
normal form} based on Garside's theory~\cite{Gar}. Here we describe 
them, successively in a non-symmetric and a symmetric version.

\subsection{Braid groups}

We first recall a few basic definitions and general results about braid groups.

\subsubsection{Presentation}

Artin's braid groups are infinite non-commu\-tative
groups. They appear in several {\it a priori} unrelated
frameworks, and they admit many equivalent definitions. 
In our case, it will be convenient to introduce them by means of explicit
presentations. 

\begin{defi}
For $n \ge 2$, the braid group~$B_n$ is defined by the
presentation
\begin{equation}\label{E:Pres}
 \langle \ss1, ..., \ss{n-1} \, ; \,
 \ss i \ss j = \ss j \ss i 
 \mbox{ for $\vert i - j\vert \ge 2$},\ 
 \ss i \ss j \ss i = \ss j \ss i \ss j 
\mbox{ for $\vert i - j\vert = 1$}
 \rangle.
\end{equation}
\end{defi}

An element of~$B_n$ will be called an {\it $n$-braid}. For each~$n$,
the identity mapping on~$\{\ss1, ..., \ss{n-1}\}$ induces an
embedding of~$B_n$ into~$B_{n+1}$, so that we can consider an
$\nn$-braid as a particular $(\nn+1)$-braid. Note that
$B_2$ is an infinite cyclic group, \ie, is isomorphic to the
group~$\Int$ of integers. For $\nn\ge3$, the group~$B_\nn$ is not
commutative: the center of~$B_n$ is the cyclic subgroup generated
by the element~$\DD\nn^2$, where $\DD\nn$ is introduced
in~\eqref{E:Delta} below.

When a group is specified using a presentation, each element of
the group is an equivalence class of words with respect to
the congruence generated by the relations of the presentation. In the sequel, a
word on the letters~$\ssss1$, ..., $\ssss{n-1}$ will be called an
{\it $n$-braid word}. So, every $n$-braid is an equivalence class of
$n$-braid words under the congruence~$\equiv$ generated
by the relations of~\eqref{E:Pres}. If the braid~$\xx$ is the
equivalence class of the word~$\ww$, we say that $\ww$ is a
representative of~$\xx$, and we write $\xx= \cl\ww$. 

\subsubsection{The word problem}

Using $\e$ for the empty word, the {\it word
problem} of~\eqref{E:Pres} is the algorithmic problem:
\begin{quote}
Given one braid word~$\ww$, does $\ww\equiv\e$ hold, \ie, does
$\ww$ represent the unit braid~$1$?
\end{quote}
This is the problem we investigate in the sequel.
Because $B_\nn$ is a group, the above one parameter problem is
equivalent to the two parameter problem:
\begin{quote}
Given two braid words~$\ww, \ww'$, does $\ww\equiv\ww'$ hold,
\ie, do $\ww$ and $\ww'$ represent the same braid?
\end{quote}
Indeed, $\ww\equiv\ww'$ is equivalent to $\ww\inv\ww'\equiv\e$,
where $\ww\inv$ is the word obtained from~$\ww$ by reversing
the order of the letters and exchanging~$\ss\ii$ and~$\sss\ii$
everywhere.

\subsubsection{Geometric interpretation}

\begin{figure}[htb]
\setlength{\unitlength}{1mm}
\begin{picture}(110,31)(0,-1)
\put(5,2){\includegraphics{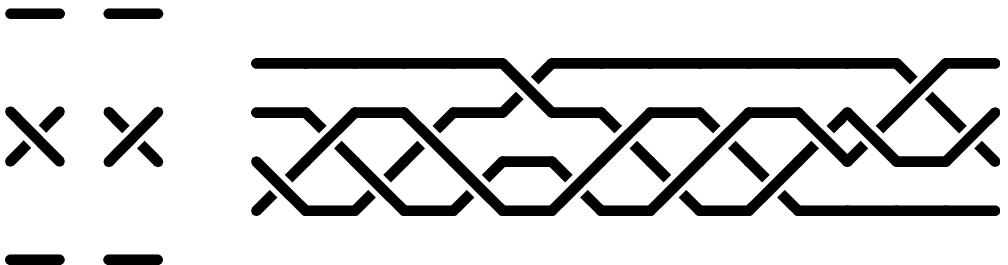}}
\put(1,1){$1$}
\put(1,11){$\ii$}
\put(-1.3,17){$\ii\!+\!1$}
\put(1,27){$\nn$}
\put(7,-2){$\ss\ii$}
\put(16,-2){$\sss\ii$}
\put(8,6){$\vdots$}
\put(18,6){$\vdots$}
\put(8,21){$\vdots$}
\put(18,21){$\vdots$}
\put(31,3){$\ss1$}
\put(35,3){$\sss2$}
\put(41,3){$\ss1$}
\put(46,3){$\ss2$}
\put(51,3){$\ss1$}
\put(56,3){$\ss3$}
\put(59.5,3){$\sss1$}
\put(65,3){$\sss2$}
\put(70.5,3){$\sss1$}
\put(76,3){$\sss2$}
\put(81.5,3){$\sss1$}
\put(87.5,3){$\ss2$}
\put(91.5,3){$\ss2$}
\put(95.5,3){$\sss3$}
\put(101,3){$\sss2$}
\end{picture}
\caption{\smaller\sf The $\nn$~strand braid diagrams
associated with~$\ss i$ and~$\sss i$, and the
$4$~strand braid diagram associated with the $4$-\nobreak braid
word
$\ss1\sss2\ss1\ss2\ss1\ss3\sss1\sss2\sss1\sss2\sss1\ss2\ss2\sss3\sss2$---also
denoted $\mathtt{aBabacABABAbbCB}$ in the sequel: one
concatenates the successive $4$~strand diagrams
corresponding to the successive letters of the word.}
\label{F:Braid}
\end{figure}

The elements of~$B_n$ can be interpreted as geometric $n$~strand
braids~\cite{Bir, KaT, Dgr}. To this end, one
associates with every braid word the plane diagram obtained by
concatenating the elementary diagrams of Figure~\ref{F:Braid}
corresponding to the successive letters. 

A braid diagram can be seen as a plane projection of a
three-dimensional figure consisting on $n$~disjoint
curves connecting the points $(1, 0, 0),
..., (n, 0, 0)$ to the points $(1, 0, 1), ..., (n, 0, 1)$ in~$\Ree^3$. Then
the relations of~\eqref{E:Pres} correspond to ambient isotopy, \ie,
to continuously moving the curves without moving their
ends and without allowing them to intersect. It is easy to check that
each relation in~\eqref{E:Pres} corresponds to such an isotopy; the
converse implication, \ie, the fact that the projections
of isotopic 3D figures can always be encoded in words connected
by~\eqref{E:Pres} was proved by E.\,Artin in~\cite{Art}. Thus the
braid word problem for the presentation~\eqref{E:Pres} is also 
the {\it braid isotopy problem}---thus similar to the (much more difficult) 
knot isotopy problem.

\subsubsection{Positive braids}

A braid word is said to be {\it positive} if it contains no
letter~$\sss\ii$. A braid is said to be {\it positive} if it can be
represented by at least one positive word. Positive $n$-braids form a
submonoid denoted~$\BB\nn$ of the group~$B_n$. Garside's theory~\cite{Gar} 
implies that $\BB\nn$ admits, as a monoid, the
presentation~\eqref{E:Pres} and that $B_n$ is a group of fractions
of~$\BB\nn$, \ie, every braid in~$B_n$ can be expressed
as~$\yy\inv\xx$ with $\xx,\yy$ in~$\BB\nn$.

For $\xx, \yy$ in~$\BB\nn$, we say that $\xx$ is a {\it
left divisor} of~$\yy$, or, equivalently,  that $\yy$ is a {\it
right multiple} of~$\xx$, denoted $\xx\divel\yy$, if $\yy = \xx
\zz$ holds for some~$\zz$ in~$\BB\nn$. Right divisors and left
multiples are defined symmetrically. With respect to left divisibility (and 
to right divisibility as well),
$\BB\nn$ has the structure of a lattice: any two positive
$\nn$-braids~$\xx, \yy$ admit a greatest common left
divisor~$\gcdl(\xx, \yy)$, and a least common right multiple. 

\subsubsection{Permutation of a braid}

The geometric interpretation makes it clear that mapping~$\ss i$
to the transposition that exchanges~$i$ and~$i+1$
induces a surjective homomorphism of the braid group~$B_n$ onto
the symmetric group~$S_n$. Under this homomorphism, here
denoted~$\perm$, a braid~$\xx$ is mapped to the permutation~$f$
of~$\{1, ..., n\}$ such that the strand that finishes at position~$i$ in
any braid diagram representing~$\xx$ begins at position~$f(i)$. 

\subsubsection{Simple braids}

A special r\^ole is played by the positive $\nn$-braid~$\DD\nn$
inductively defined by
\begin{equation} 
\label{E:Delta}
\DD1 = 1, \qquad
\DD\nn = \ss1 \ss2 \dots \ss{\nn-1} \, \DD{\nn-1}.
\end{equation}
In~$\BB\nn$, the left and the right divisors of~$\DD\nn$ coincide,
they include each of~$\ss1$, ..., $\ss{\nn-1}$, and they make
a finite sublattice of~$\BB\nn$ with $\nn!$~elements. These divisors
of~$\DD\nn$ are called {\it simple} braids. Geometrically, simple
braids are those positive braids that can be represented by a braid
diagram in which any two strands cross at most once. Moreover, the
restriction of the projection~$\perm$ to simple braids is a bijection:
for each permutation~$\ff$ in~$S_\nn$, there exists exactly
one simple braid~$\xs$ satisfying $\perm(\xs)=\ff$. This simple
braid will be denoted by~$\br\ff$.

\subsection{The greedy normal form}
\label{S:Greedy}

The seminal results of F.A.\,Garside~\cite{Gar} subsequently
developed in~\cite{Thu, Adj, ElM} imply that braid groups can be
equipped with a remarkable normal form, the so-called
greedy normal form. The latter  is excellent both in theory and in
practice as it provides a  bi-automatic structure, and it is easily
computable. 

\subsubsection{Description}

The greedy normal form exists in several variants, in particular left
and symmetric right versions. Here we shall consider the left versions
only. Again, there exist two different versions. Both consist in
expressing an arbitrary braid as a quotient of two positive braids,
\ie, as a fraction. In the version considered in this section, all
denominators have some special form, namely they are powers of
the element~$\DD\nn$. By contrast, in the version
considered in Section~\ref{S:Symmetric} below, the numerator and
the denominator of fractions play symmetric r\^oles.

\begin{defi}
\label{D:Normal}
$(i)$ A sequence of simple braids $(\xs_1, ..., \xs_\pp)$ is said to be
{\it normal} if, for each~$\kk<\pp$, every~$\ss\ii$ that is a left
divisor of~$\xs_{\kk+1}$ is a right divisor of~$\xs_\kk$.

$(ii)$ A sequence of permutations $(\ff_1, ..., \ff_\pp)$ is said to be
{\it normal} if, for each~$\kk<\pp$,  every recoil
of~$\ff_{\kk+1}$ is a descent of~$\ff_\kk$, \ie, if
$\ff_{\kk+1}\inv(i) > \ff_{\kk+1}\inv(i+1)$ implies $\ff_\kk(i) >
\ff_\kk(i+1)$.
\end{defi}

The connection between~$(i)$ and~$(ii)$ in
Definition~\ref{D:Normal} is that, if $\xs$ is a simple $\nn$-braid
and $\ff$ is the associated permutation, then $\ss\ii$ is a left (\resp
right) divisor of~$\xs$ if and only if we
have $\ff\inv(\ii)>\ff\inv(\ii+1)$ (\resp $\ff(\ii)> \ff(\ii+1)$).
So a sequence of simple braids $(\xs_1, ..., \xs_\pp)$ is normal if and
only if the associated sequence of permutations $(\perm(\xs_1),
...,\perm(\xs_\pp))$ is normal.

We denote by~$\omega_\nn$ the flip permutation of~$\{1, ...,
\nn\}$ defined by $\omega_\nn(\ii) = \nn-\ii+1$.

\begin{theo}
\cite[Chapter~9]{Eps}
\label{T:Normal}
$(i)$ Every braid~$\zz$ in~$B_n$ admits a unique decomposition of
the form $\DD\nn^\mm \xs_1 ... \xs_\pp$ with $\mm$ in~$\Int$
and $(\xs_1, ..., \xs_\pp)$ a normal sequence of simple braids
satisfying $\xs_1 \not=\DD\nn$ and $\xs_\pp\not=1$.

$(ii)$ Every braid~$\zz$ in~$B_n$ admits a unique decomposition of
the form $\DD\nn^\mm \br{\ff_1} ... \br{\ff_\pp}$ with $\mm$
in~$\Int$ and $(\ff_1, ..., \ff_\pp)$ a normal sequence of
permutations satisfying $\ff_1 \not=\omega_\nn$ and
$\ff_\pp\not=\id$.
\end{theo}

In the situation of Theorem~\ref{T:Normal}$(i)$, the sequence
$(\mm; \xs_1, ..., \xs_\pp)$ is called the {\it greedy normal form}
of~$\zz$---or the $\nn$-greedy normal form of~$\zz$ if we wish to insist
on the braid index~$\nn$. As simple braids are in
one-to-one correspondence with permutations, and by the remark
above, the braid form and the permutation form of the greedy normal
form are equivalent. So there is no problem
in also calling the sequence
$(\mm; \ff_1, ..., \ff_\pp)$ of Theorem~\ref{T:Normal}$(ii)$ the
greedy normal form of~$\zz$. 

Clearly, $(0; \emptyset)$ is the greedy normal form of~$1$, and the
uniqueness of the greedy normal form implies the following solution to the
braid isotopy problem---but a solution that remains uneffective as long 
as we give no method for computing from an arbitrary braid word~$\ww$
the greedy normal form of~$\cl\ww$, \ie,
until Section~\ref{S:Grid} below:

\begin{coro}
\label{C:Greedy}
A braid word~$\ww$ represents~$1$ in the braid group if and
only if the greedy normal form of~$\cl\ww$ is~$(0; \emptyset)$.
Two braid words~$\ww, \ww'$ represent the same braid
in~$B_\nn$ if and only if the greedy normal forms of the
braids~$\cl\ww$ and~$\cl{\ww'}$ coincide.
\end{coro}

\begin{exam}
\label{X:NF}
In order to obtain shorter notation, we shall in the sequel use $\tta,
\ttb, \ttc...$ for~$\ss1, \ss2, \ss3...$, and, symmetrically,
$\ttA, \ttB...$ for $\sss1, \sss2$... (as in the caption of
Figure~\ref{F:Braid}). Then, a typical greedy normal form for a
$4$-braid is the sequence

$(-2; \mathtt{ac, abcb, bcba, a})$,

\noindent\ie, equivalently, using $(\ff(1), ..., \ff(\nn))$ to
specify a permutation~$\ff$ of~$\{1, ..., \nn\}$,

$(-2; (2,1,4,3), (2,4,3,1), (4,1,3,2), (2,1,3,4))$,

\noindent consisting of an integer and four simple $4$-braids, or,
equivalently, four permutations of~$\{1, ..., 4\}$: for instance, 
$(2,1,4,3)$ is the permutation associated
with~$\mathtt{ac}$, \ie, with~$\ss1\ss3$. To check that we have a greedy 
normal form, we observe for instance that the descents of  $(2,1,4,3)$ 
are $1$ and~$3$, while the recoils of $(2,4,3,1)$, \ie, the descents of~$(4,1,2,3)$, are
$1$ and~$3$ as well, so the normality condition is satisfied between
$(2,1,4,3)$ and~$(2,4,3,1)$. The other verifications are similar.
So, the above sequences are two versions of the greedy normal
form of the $4$-braid represented by $\Delta_4^{-2}.\mathtt{ac.abcb.bcba.a}$, \ie,
by

$\ww = \mathtt{ABACBA.ABACBA.ac.abcb.bcba.a}$.

\noindent As the above normal form is not $(0, \emptyset)$, we
deduce from Corollary~\ref{C:Greedy} that $\ww$ does not represent~$1$ in~$B_4$.
\end{exam}

\subsubsection{Explanation}

The existence and uniqueness of the greedy normal form follows
from two results: 

$(i)$ For every braid~$\zz$ in~$B_\nn$, there exist $\mm\in\Int$ 
and $\xx\in\BB\nn$ satisfying $\DD\nn^\mm\xx=\zz$, the
decomposition being unique provided $\DD\nn\not\divel\xx$ is
required;

$(ii)$ For every positive braid~$\xx$, there exists a unique normal
sequence $(\xs_1, ..., \xs_\pp)$ of simple $\nn$-braids with
$\xs_\pp\not=1$ satisfying $\xs_1 ... \xs_\pp = \xx$.

The proof of~$(i)$ is an easy induction on the length of a braid
word representing~$\zz$ once one knows, for each~$\ii$, the
relations $\ss\ii\divel\DD\nn$ and $\DD\nn
\ss\ii=\ss{\nn-\ii}\DD\nn$, which imply that, by multiplying an
$\nn$-braid by a sufficient large power of~$\DD\nn$, one can
always obtain a positive braid.

As for~$(ii)$, the existence of left gcd's in the monoid~$\BB\nn$
implies that each positive $\nn$-braid~$\xx$ can be expressed as
$\xx = \xs_1\xx'$ with $\xs_1 = \gcdl(\xx, \DD\nn)$, and
$\xx\not=1$ implies $\xs_1\not=1$. By iterating the process,
thus writing $\xx' = \xs_2\xx''$, {\it etc.}, one eventually obtains a
decomposition $\xx= \xs_1 ... \xs_\pp$. By construction, the
sequence $(\xs_1, ..., \xs_\pp)$ consists of divisors of~$\DD\nn$,
\ie, of simple $\nn$-braids, and, for each~$\kk<\pp$, one has
$\xs_\kk = \gcdl(\xs_\kk
\xs_{\kk+1} ... \xs_\pp, \DD\nn)$, hence, {\it a fortiori}, $\xs_\kk =
\gcdl(\xs_\kk \xs_{\kk+1}, \DD\nn)$. The point is that the
latter relations are equivalent to $(\xs_\kk, \xs_{\kk+1})$
being normal for each~$\kk$, and, therefore, to $(\xs_1, ...,
\xs_\pp)$ being normal. The uniqueness comes from the fact that, if
$(\xs_1, ..., \xs_\pp)$ is a normal sequence, then, necessarily, one
has $\xs_1 = \gcdl(\xs_1...\xs_\pp, \DD\nn)$.

\subsubsection{Discussion}

What is missing in the above description of the greedy
normal form is an algorithm for computing the (unique)
normal form of a braid~$\zz$ starting from an arbitrary
representative of~$\zz$. Clearly, the existence of such an algorithm is
a necessary condition for using the normal form in
practice. Such an algorithm will be provided in
Section~\ref{S:Grid} below, and discussing the practical
implementation of the greedy normal form will be possible only then.

Actually, the point is not necessarily to find the normal form
equivalent to an arbitrary word, but rather to find the normal form of
the product or the quotient of two normal forms. Indeed, whenever
one chooses to work with normal forms, one may forget about
non-normal words provided one is able to perform the basic
operations inside the family of normal forms. Of course, as the
generator~$\ss\ii$ is itself a braid, with normal form~$(0; \ss\ii)$,
an algorithm computing the product of two normal forms will in
particular determine the product of a normal form by~$\ss\ii$ and,
therefore, inductively determine the normal form of any product
of~$\ssss\ii$'s, but the general philosophy is not exactly that of a
normalizing algorithm.

Note that, while the permutation variant of the greedy normal form is 
non-ambiguously defined, the
braid word variant is not: for instance, the first braid factor in 
Example~\ref{X:NF}, namely~$\mathtt{ac}$, is uniquely
defined as a simple braid, but it can be represented by two different
positive braid words, namely $\mathtt{ac}$ and~$\mathtt{ca}$. So
the braid word form becomes unique only when a distinguished
word representative has been chosen for every simple braid. This
explains why the permutation form is often more convenient.

\subsection{The symmetric normal form}
\label{S:Symmetric}

In the greedy normal form where the denominator is always a power of~$\DD\nn$. 
The symmetric normal form is a variant in which the numerator and the 
denominator play symmetric r\^oles.

\subsubsection{Description}

The symmetric normal form appeals to the same notion of a normal
sequence of simple braids as the greedy normal form.

\begin{theo}
\cite[Chapter~9]{Eps}
\label{T:Symmetric}
$(i)$ Every braid~$\zz$ admits a unique decomposition as
$\xt_\qq\inv\, ... \,\xt_1\inv \,\xs_1\, ... \,\xs_\pp$ with $(\xs_1, ...,
\xs_\pp)$, $(\xt_1, ..., \xt_\qq)$ two normal sequences of simple braids
satisfying $\xs_\pp\not=1$, $\xt_\qq\not=1$, and
$\gcdl(\xs_1, \xt_1)=1$.

$(ii)$ Every braid~$\zz$ admits a unique decomposition
as $\br{\gg_\qq}{}\inv \,...\, \br{\gg_1}{}\inv\,\br{\ff_1}\, ...
\,\br{\ff_\pp}$ with $(\ff_1, ..., \ff_\pp)$, $(\gg_1, ..., \gg_\qq)$ two
normal sequences of permutations satisfying
$\ff_\pp, \gg_\qq \not= \id$ and
such that $\ff_1\inv(\ii)>\ff_1\inv(\ii+1)$ implies
$\gg_1\inv(\ii)\le\gg_1\inv(\ii+1)$.
\end{theo}

In the situation of Theorem~\ref{T:Symmetric}$(i)$, the
double sequence $(\xt_1, ..., \xt_\qq; \xs_1, ..., \xs_\pp)$ is called
the {\it symmetric normal form} of~$\zz$. As for the greedy normal
form, both versions of the symmetric normal forms are equivalent,
and the double sequence of permutations $(\gg_1, ..., \gg_\qq; \ff_1,
..., \ff_\pp)$ of Theorem~\ref{T:Symmetric}$(ii)$ is also called the
symmetric normal form of~$\zz$. As $(\emptyset;
\emptyset)$ is the symmetric normal form of~$1$, the uniqueness 
of the symmetric normal form provides the following
solution to the braid isotopy problem:

\begin{coro}
\label{C:Symmetric}
A braid word~$\ww$ represents~$1$ in the braid group if and
only if the symmetric normal form of~$\cl\ww$ is~$(\emptyset;
\emptyset)$. Two braid words~$\ww, \ww'$ represent the same braid
in~$B_\nn$ if and only if the symmetric normal forms of~$\cl\ww$
and~$\cl{\ww'}$ coincide.
\end{coro}

\begin{exam}
\label{X:Symmetric}
A typical symmetric normal form is

$(\mathtt{ab, bacb; bcba, a})$,

\noindent\ie, equivalently, 

$((2,3,1,4), (3,4,1,2); (4,1,3,2), (2,1,3,4))$,

\noindent For instance, the normality condition between the first factors of the
two sequences holds as $1$ is the only recoil of~$(2,3,1,4)$, while
$2$ and~$3$ are those of~$(4,1,3,2)$. In other words, the simple
braids~$\mathtt{ab}$ and~$\mathtt{bcba}$ admit no nontrivial
common left divisor. Thus the above expressions specify the
symmetric normal form of the braid represented by

$\ww = \mathtt{BCAB.BA.bcba.a}$,

\noindent which we shall see below coincides with the braid of
Example~\ref{X:NF}---and of Figure~\ref{F:Braid}. As the above
normal form is not $(\emptyset, \emptyset)$, we deduce that $\ww$
 does not represent~$1$ in~$B_4$.
\end{exam}

\subsubsection{Explanation}

As $B_\nn$ is a group of fractions for the monoid~$\BB\nn$, every
element of~$B_\nn$ can be expressed as a fraction~$\yy\inv\xx$
with $\xx, \yy$ in~$\BB\nn$, and the decomposition is unique if, in
addition, $\gcdl(\xx,\yy)$ is required. The symmetric normal form
is obtained by taking the greedy normal forms of the positive
braids~$\xx, \yy$ so obtained, with the only difference that the 
$\DD\nn$ factors are not separated. The specific properties of 
normal sequences imply that $\gcdl(\xx,\yy)=1$ is
equivalent to $\gcdl(\xs_1, \xt_1)=1$, where $\xs_1$ and $\xt_1$
respectively are the first factors in the normal forms of~$\xx$
and~$\yy$.

\subsubsection{Discussion}

As in Section~\ref{S:Greedy}, our description will be complete only
when we give algorithms for computing the symmetric normal form
of a product or of a quotient. This will be done in
Section~\ref{S:Grid}.

\section{Grid properties of the normal form}
\label{S:Grid}

The interest of the greedy and symmetric normal forms of braids
lies in the existence of simple computing rules for determining the
normal form of a product or of a quotient. Here we shall solve the
following two problems, which then easily leads to
complete algorithms for all problems connected with the
normal forms:
\begin{quote}
Starting with the greedy normal forms of two positive
braids~$\xx, \yy$, find the greedy normal form of the
product~$\yy\xx$ and, assuming that $\yy$ is a right divisor
of~$\xx$, the greedy normal form of the quotient~$\xx\yy\inv$.
\end{quote}
The solutions we shall describe involve grid diagrams that visualize
the properties of normal sequences. To this end, it will be convenient
to associate with every braid~$\xx$ an arrow labelled~$\xx$, so
that a  relation of the form $\xx\yy=\zz$ corresponds to a
commutative diagram for the associated arrows. We
shall indicate that a sequence $(\xs_1, ...,
\xs_\pp)$ is normal by drawing an  arc connecting the final end of
each arrow with the initial end of the next one, on the shape of
\begin{center}
\begin{picture}(75,6)
\put(0,0){\includegraphics{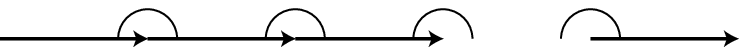}}
\put(6,2){$\xs_1$}
\put(21,2){$\xs_2$}
\put(36,2){$\xs_3$}
\put(50,0.5){$\ldots$}
\put(66,2){$\xs_\pp$}
\end{picture}.
\end{center}

\subsection{Prerequisites}

The only properties of simple braids needed in the forthcoming
proofs are summarized in the following lemma, which is standard:

\begin{lemm}
\cite{Cha}
\label{L:Head}
For~$\xx$ in~$\BB\nn$, let $\HH(\xx):=\gcdl(\xx, \DD\nn)$. Then,
for all~$\xx, \yy$, we have:
\begin{gather}
\label{E:Head2}
\HH(\xx) \divel \xx,
\text{\quad and $\HH(\xx)=\xx$ if $\xx$ is simple},\\
\label{E:Head3}
\xx\divel\yy  \text{ implies } 
\HH(\xx) \divel \HH(\yy),\\
\label{E:Head4}
\HH(\xx\yy) =\HH(\xx\HH(\yy)).
\end{gather}
\end{lemm}

Then the definition of normal sequences immediately rewrites in
terms of the function~$\HH$ as follows---and this is the
technical form we shall use in the sequel:

\begin{lemm}
A sequence of simple braid~$(\xs_1, ..., \xs_\pp)$ is normal if and
only if, for each~$\kk<\pp$, we have $\xs_\kk = \HH(\xs_\kk
\xs_{\kk+1})$.
\end{lemm}

In the sequel, we mostly deal with positive braids. When we say that
$(\xs_1, ..., \xs_\pp)$ is the normal form of a positive $\nn$-braid~$\xx$, we
mean that $(\xs_1, ..., \xs_\pp)$ is normal and that $\xx = \xs_1 ... \xs_\pp$.
Equivalently, this means that the greedy normal form of~$\xx$ is
$(\mm; \xs_{\mm+1}, ..., \xs_\pp)$, where $\mm$ is the number of initial~$\kk$'s
such that $\xs_\kk$ equals~$\DD\nn$.

\subsection{Grid properties for the quotient}
\label{S:GridComplement}

We start with the problem of finding the normal form
of~$\xx\yy\inv$ from those of~$\xx$ and~$\yy$, where $\xx$
and~$\yy$ belong to~$\BB\nn$. In principle, the computation
makes sense only if $\yy$ happens to be a right divisor of~$\xx$.
Actually, in any case, the positive braids~$\xx, \yy$ admit a left lcm
in~$\BB\nn$, and what we shall do is to compute the normal form of
this left lcm, denoted~$\lcml(\xx,\yy)$, and of the associated left 
complements defined as follows:

\begin{defi}
For $\xx,\yy$ in~$\BB\nn$, the unique~$\xx'$ in~$\BB\nn$
satisfying $\lcml(\xx,\yy)=\xx'\yy$ is called the {\it left
complement} of~$\xx$ in~$\yy$, and denoted by~$\xx/\yy$.
\end{defi}

If $\yy$ happens to be a right divisor of~$\yy$, \ie, if
$\yy/\xx = 1$ holds, then we have $\lcml(\xx,\yy) = \yy$,
and $\xx/\yy = \xx\yy\inv$. Thus an algorithm
computing the left complement is in particular an algorithm
computing the right quotient when it exists.

A standard observation is that simple braids are closed under left
lcm and left complement. Indeed, assume that $\xs,
\xt$ are simple, and let $\xu = \lcml(\xs, \xt) = \xs' \xt = \xt' \xs$,
\ie, $\xs'= \xs/\xt$ and $\xt' = \xt/\xs$. Then, $\xs$ and $\xt$ are
right divisors of~$\DD\nn$, hence $\xu$ is also a right divisor
of~$\DD\nn$, and it is therefore a simple braid. Then, $\xs'$
and~$\xt'$ are simple as well, as every left divisor of a
simple braid is a simple braid. Moreover, it is straightforward to
check that $\xs'\xt$ being the left lcm of~$\xs$ and~$\xt$ is
equivalent to $\xs'$ and $\xt'$ admitting no common left divisor
except~$1$, \ie, to $\gcdl(\xs', \xt')=1$.

\begin{figure}[htb]
\begin{picture}(23,19)
\put(0,.5){\includegraphics{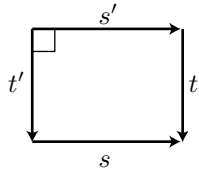}}
\put(-2,8){$\xt'$}
\put(22,8){$\xt$}
\put(10,-2){$\xs$}
\put(10,17){$\xs'$}
\end{picture}
\caption{\smaller\sf A $C$-tile: a commutative diagram involving
four simple braids such that the two initial ones have no nontrivial
common left divisor---here indicated by a perpendicularity sign.}
\label{F:CTile}
\end{figure}

\begin{defi}
A commutative diagram consisting of four simple braids~$\xs, \xt,
\xs', \xt'$ satisfying $\xs'\xt=\xt'\xs$ and $\gcdl(\xs', \xt')=1$ is
called a {\it $C$-tile} (like ``complement tile''), and it is represented
as in Figure~\ref{F:CTile}.
\end{defi}

We establish normality conditions for diagrams involving
the above $C$-tiles.

\begin{figure}[htb]
\begin{picture}(41,24)
\put(0,1){\includegraphics{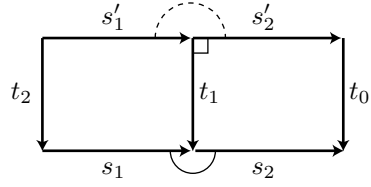}}
\put(-3,11){$\xt_2$}
\put(22,11){$\xt_1$}
\put(42,11){$\xt_0$}
\put(9,1){$\xs_1$}
\put(29,1){$\xs_2$}
\put(9,21){$\xs'_1$}
\put(29,21){$\xs'_2$}
\end{picture}
\caption{\smaller\sf A grid property involving $C$-tiles:
assuming that the right rectangle is a $C$-tile and the bottom
sequence $(\xs_1, \xs_2)$ is normal, the top sequence $(\xs'_1,
\xs'_2)$ is normal as well.}
\label{F:Complement1}
\end{figure}

\begin{lemm}
\label{L:Complement1}
(Figure~\ref{F:Complement1}) 
Assume that $\xs_1, \xs_2, \xs'_1, \xs'_2, \xt_0, \xt_1,
\xt_2$ are simple braids satisfying $\xt_2\xs_1=\xs'_1\xt_1$, 
$\xt_1\xs_2=\xs'_2\xt_0$, $\gcdl(\xs'_2,\xt_1)=1$, and that
$(\xs_1,\xs_2)$ is normal. Then $(\xs'_1, \xs'_2)$ is normal as well.
\end{lemm}

\begin{proof}
With the notation of Figure~\ref{F:Complement1}, we compute:
\begin{align*}
\HH(\xs'_1\xs'_2)
&\divel \HH(\xs'_1\xs'_2\xt_0)
&&\text{by~\eqref{E:Head3}}\\
&= \HH(\xt_2\xs_1\xs_2)
&&\text{by commutativity}\\
&= \HH(\xt_2\HH(\xs_1\xs_2))
&&\text{by~\eqref{E:Head4}}\\
&= \HH(\xt_2 \xs_1)
&&\text{by the hypothesis that $(\xs_1,\xs_2)$ is normal}\\
&= \HH(\xs'_1\xt_1)
&&\text{by commutativity}\\
&\divel \xs'_1\xt_1
&&\text{by~\eqref{E:Head2}}.
\end{align*}
On the other hand, by~\eqref{E:Head2} again, we have
$\HH(\xs'_1\xs'_2) \divel \xs'_1\xs'_2$, hence
$$\HH(\xs'_1\xs'_2) \divel \gcdl(\xs'_1\xs'_2, \xs'_1\xt_1) =
\xs'_1 \cdot \gcdl(\xs'_2, \xt_1) = \xs'_1,$$
using the hypothesis $\gcdl(\xs'_2, \xt_1)=1$. It follows that 
$(\xs'_1, \xs'_2)$ is normal.
\end{proof}

Similarly, we have the following result involving the diagonals of
$C$-tiles, \ie, the left lcm's of the corresponding simple braids.

\begin{figure}[htb]
\begin{picture}(41,35)
\put(0,1){\includegraphics{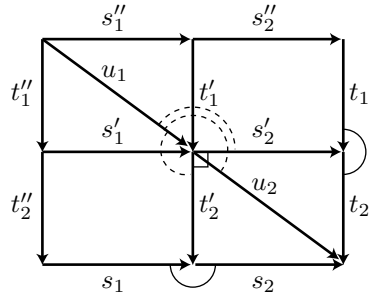}}
\put(-3,11){$\xt''_2$}
\put(22,11){$\xt'_2$}
\put(42,11){$\xt_2$}
\put(-3,26){$\xt''_1$}
\put(22,26){$\xt'_1$}
\put(42,26){$\xt_1$}
\put(9,1){$\xs_1$}
\put(29,1){$\xs_2$}
\put(9,21){$\xs'_1$}
\put(29,21){$\xs'_2$}
\put(9,36){$\xs''_1$}
\put(29,36){$\xs''_2$}
\put(9,29){$\xu_1$}
\put(29,14){$\xu_2$}
\end{picture}
\caption{\smaller\sf Another grid property involving $C$-tiles:
assuming that the bottom right rectangle is a $C$-tile and the
bottom and right sequences $(\xs_1, \xs_2)$ and $(\xt_1, \xt_2)$
are normal, the diagonal sequence $(\xu_1, \xu_2)$ is normal
as well, and so are the diagonal-then-horizontal sequence
$(\xu_1,\xs'_2)$ and the diagonal-then-vertical sequence
$(\xu_1,\xt'_2)$.}
\label{F:Complement2}
\end{figure}

\begin{lemm}
\label{L:Complement2}
(Figure~\ref{F:Complement2})
Assume that $(\xs_1, \xs_2)$ and $(\xt_1, \xt_2)$ are normal
sequences of simple braids. Let $\xu_1, \xu_2$ be as in
Figure~\ref{F:Complement2}. Then $(\xu_1, \xu_2)$, $(\xu_1,
\xs'_2)$, and $(\xu_1, \xt'_2)$ are normal as well. 
\end{lemm}

\begin{proof}
With the notation of Figure~\ref{F:Complement2}, we compute
\begin{align*}
\HH(\xu_1\xu_2)
&= \HH(\xt''_2\xt''_1\xs_1\xs_2)
&&\text{by commutativity}\\
&= \HH(\xt''_2\xt''_1\HH(\xs_1\xs_2))
&&\text{by~\eqref{E:Head4}}\\
&= \HH(\xt''_2\xt''_1\xs_1)
&&\text{by the hypothesis that $(\xs_1, \xs_2)$ is normal}\\
&= \HH(\xu_1\xt'_2)
&&\text{by commutativity}\\
&\divel \xu_1\xt'_2
&&\text{by~\eqref{E:Head2}}.
\end{align*}
Similarly, we obtain $\HH(\xu_1\xu_2)\divel \xu_1\xs'_2$,
and we deduce
$$\HH(\xu_1\xu_2)\divel \gcdl(\xu_1\xt'_2,
\xu_1\xs'_2) = \xu_1 \cdot \gcdl(\xs'_2, \xt'_2) = \xu_1$$ from
the hypothesis $\gcdl(\xs'_2, \xt'_2)=1$. So 
$(\xu_1,\xu_2)$ is normal. As $\xs'_2$ is a left divisor
of~$\xu_2$, the normality of~$(\xu_1, \xu_2)$ implies
that of~$(\xu_1, \xs'_2)$, and, similarly, that of~$(\xu_2,
\xt'_2)$.
\end{proof}

\begin{figure}
\begin{picture}(110,45)
\put(0,-0.5){\includegraphics{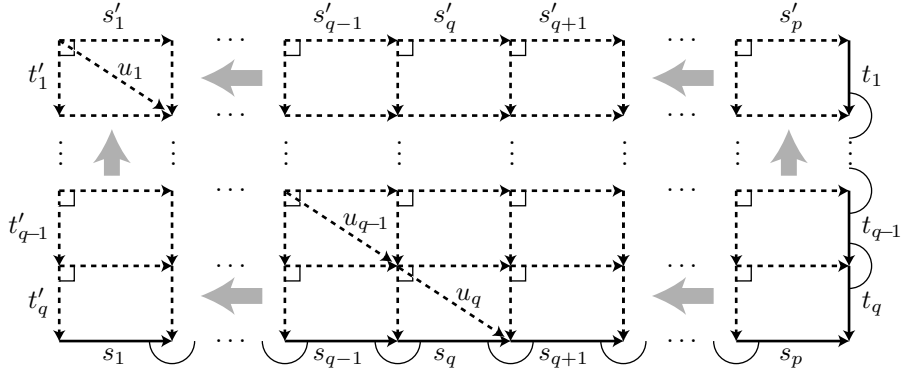}}
\put(8,0){$\xs_1$}
\put(36,0){$\xs_{\qq-1}$}
\put(52,0){$\xs_{\qq}$}
\put(66,0){$\xs_{\qq+1}$}
\put(98,0){$\xs_{\pp}$}
\put(8,45){$\xs'_1$}
\put(36,45){$\xs'_{\qq-1}$}
\put(52,45){$\xs'_{\qq}$}
\put(66,45){$\xs'_{\qq+1}$}
\put(98,45){$\xs'_{\pp}$}
\put(1,7){\llap{$\xt'_\qq$}}
\put(1,17){\llap{$\xt'_{\qq\!-\!1}$}}
\put(1,37){\llap{$\xt'_1$}}
\put(109,7){$\xt_\qq$}
\put(109,17){$\xt_{\qq\!-\!1}$}
\put(109,37){$\xt_1$}
\put(2,26){$\vdots$}
\put(17,26){$\vdots$}
\put(32,26){$\vdots$}
\put(47,26){$\vdots$}
\put(62,26){$\vdots$}
\put(77,26){$\vdots$}
\put(92,26){$\vdots$}
\put(107,26){$\vdots$}
\put(23,2.5){$\ldots$}
\put(23,12.5){$\ldots$}
\put(23,22.5){$\ldots$}
\put(23,32.5){$\ldots$}
\put(23,42.5){$\ldots$}
\put(83,2.5){$\ldots$}
\put(83,12.5){$\ldots$}
\put(83,22.5){$\ldots$}
\put(83,32.5){$\ldots$}
\put(83,42.5){$\ldots$}
\put(10,38){$\xu_1$}
\put(40,18){$\xu_{\qq\!-\!1}$}
\put(55,8){$\xu_{\qq}$}
\end{picture}
\caption{\smaller\sf Construction of the left lcm and the
left complements by the grid method: we start from the
bottom and right sides, and fill the grid with $C$-tiles, from right
to left and from bottom to top; then every row, every column, and
every diagonal in the grid consist of normal sequences.}
\label{F:Complement}
\end{figure}

We are ready to establish:

\begin{prop}
\label{P:Complement}
$(i)$ Assume that $(\xs_1, ..., \xs_\pp)$ and $(\xt_1, ..., \xt_\qq)$
are normal sequences of simple braids. Let $\mathcal{D}$ be the grid
diagram obtained by starting from the right column $(\xt_1, ...,
\xt_\qq)$ and the bottom row $(\xs_1, ..., \xs_\pp)$ and filling the
diagram with $C$-tiles from right to left, and from bottom to top,
 as shown in Figure~\ref{F:Complement}.
Then every path in~$\mathcal{D}$ consisting of diagonal arrows
followed by horizontal arrows, as well as every path consisting of
diagonal arrows followed by horizontal arrows  corresponds to a
normal sequence. 

$(ii)$ Let $\xx = \xs_1...\xs_\pp$ and $\yy=\xt_1...\xt_\qq$. Let
$(\xs'_1, ..., \xs'_\pp)$ be the top row of~$\mathcal{D}$, 
let $(\xt'_1, ..., \xt'_\qq)$ be its left column, and, assuming
$\pp\ge\qq$, let $(\xu_1, ..., \xu_\qq)$ be the diagonal from the
top-left corner. Then $(\xs'_1, ..., \xs'_\pp)$ is the normal form
of~$\xx/\yy$, while $(\xt'_1, ..., \xt'_\qq)$ is the normal form
of~$\yy/\xx$, and $(\xu_1, ..., \xu_\qq, \xs_{\qq+1}, ..., \xs_\pp)$
is the normal form of the left lcm of~$\xx$ and~$\yy$.
\end{prop}

\begin{proof}
By construction, the diagram~$\mathcal{D}$ is commutative. Let
$\zz = \xs'_1 ... \xs'_\pp \xt_1 ... \xt_\qq$. Then $\zz$ is a
common left multiple of~$\xx$ and~$\yy$. Moreover, an easy
induction shows that every common left multiple~$\zz'$ of~$\xx$
and~$\yy$ has to be a left multiple of~$\zz$: indeed, $\zz'$, being
a common left multiple of~$\xs_\pp$ and~$\xt_\qq$, has to be a
left multiple of the braid represented by the diagonal of the
bottom-right rectangle in~$\mathcal{D}$, and we argue similarly for
each of the $\pp\qq$~rectangles in~$\mathcal{D}$. Hence $\zz$ is
the left lcm of~$\xx$ and~$\yy$, and, therefore, we have $\xs'_1 ...
\xs'_\pp = \xx/\yy$ and $\xt'_1 ... \xt'_\qq = \yy/\xx$.

Then, by repeatedly applying Lemma~\ref{L:Complement1}, we
obtain that every row and every column in~$\mathcal{D}$ is a
normal sequence. As for the diagonals, and the diagonals followed by
rows or columns, we similarly apply Lemma~\ref{L:Complement2}.
\end{proof}

As a straightforward application, we obtain the following computing
rule for determining the normal form of~$\xx\xu\inv$ from that
of~$\xx$.

\begin{figure}[htb]
\begin{picture}(100,22)
\put(0,0.7){\includegraphics{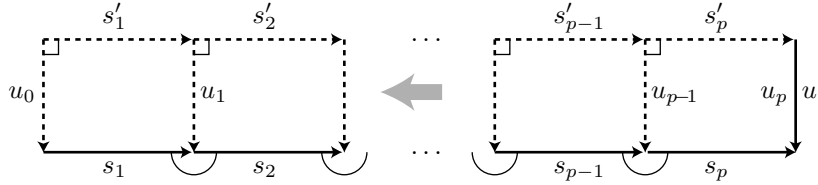}}
\put(9,1){$\xs_1$}
\put(29,1){$\xs_2$}
\put(69,1){$\xs_{\pp-1}$}
\put(89,1){$\xs_\pp$}
\put(-3.5,11){$\xu_0$}
\put(22,11){$\xu_1$}
\put(82,11){$\xu_{\pp\!-\!1}$}
\put(96.5,11){$\xu_\pp$}
\put(102,11){$\xu$}
\put(9,21){$\xs'_1$}
\put(29,21){$\xs'_2$}
\put(69,21){$\xs'_{\pp-1}$}
\put(89,21){$\xs'_\pp$}
\put(50,3.5){$\ldots$}
\put(50,18.5){$\ldots$}
\end{picture}
\caption{\smaller\sf Normal form of~$\xx\xu\inv$ from
that of~$\xx$: start from the plain arrows and fill the
diagram with $C$-tiles, from right to left; the expected normal
form is the sequence $(\xs'_1, ..., \xs'_\pp)$ read on the top line; the
hypothesis that $\xu$ is a right divisor of~$\xx$ guarantees that
$\xu_0$ is~$1$---without any hypothesis, $(\xs'_1, ..., \xs'_\pp)$ is
the normal form of~$\xx/\xu$, and $\xu_0$ is then~$\xu/\xx$.}
\label{F:ComplementOne}
\end{figure}

\begin{coro}
\label{C:Complement}
Assume that $(\xs_1, ..., \xs_\pp)$ is the normal form of
a positive braid~$\xx$, and that $\xu$ is a simple braid that
divides~$\xx$ on the right. Then the normal form of~$\xx\xu\inv$
is the sequence $(\xs'_1, ..., \xs'_\pp)$ determined by
the grid diagram of Figure~\ref{F:ComplementOne}.
\end{coro}

\subsection{Grid properties for the product}
\label{S:GridProduct}

We turn to the product, with the question of determining the normal
form of~$\yy\xx$ from the normal forms of~$\xx$ and~$\yy$. The
method is similar to that of the previous section, with another type
of basic tile.

For all simple braids~$\xt_1, \xt_2$, there exist simple
braids~$\xs_1, \xs_2$ satisfying $\xs_1 \xs_2 = \xt_1 \xt_2$ and
$(\xs_1, \xs_2)$ is normal. Indeed, let $\xs_1=\HH(\xt_1\xt_2)$,
and~$\xs_2$ satisfying $\xt_1\xt_2 = \xs_1\xs_2$.
By~\eqref{E:Head2}, $\xt_1\divel\xt_1\xt_2$ implies $\xt_1 =
\HH(\xt_1) \divel \HH(\xt_1\xt_2)=\xs_1$, so we have
$\xt_2=\xu\xs_2$ for some~$\xu$,  which forces~$\xs_2$ to be
simple. In other words, the normal form of the product of two
simple braids must consist of at most two simple braids.

\begin{figure}[htb]
\begin{picture}(21,16)
\put(0,0.2){\includegraphics{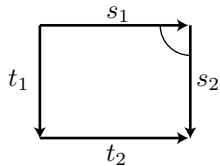}}
\put(-3,8){$\xt_1$}
\put(22,8){$\xs_2$}
\put(10,-2){$\xt_2$}
\put(10,17){$\xs_1$}
\end{picture}
\caption{\smaller\sf A $P$-tile: a commutative diagram involving
four simple braids such that two of them make a normal sequence.}
\label{F:PTile}
\end{figure}

\begin{defi}
A commutative diagram consisting of four simple braids~$\xs_1$,
$\xs_2$, $\xt_1$, $\xt_2$ satisfying $\xs_1 \xs_2 = \xt_1 \xt_2$
and such that $(\xs_1, \xs_2)$ is normal is called a {\it $P$-tile}
(like ``product tile''), and it is represented as in
Figure~\ref{F:PTile}.
\end{defi}

As was done above with $C$-tiles, we establish normality results
in diagrams involving $P$-tiles.

\begin{figure}[htb]
\begin{picture}(41,21)
\put(0,0){\includegraphics{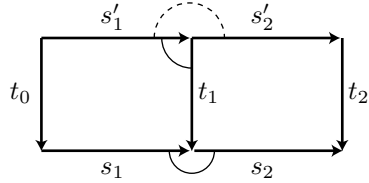}}
\put(-3,10){$\xt_0$}
\put(22,10){$\xt_1$}
\put(42,10){$\xt_2$}
\put(9,0){$\xs_1$}
\put(29,0){$\xs_2$}
\put(9,20){$\xs'_1$}
\put(29,20){$\xs'_2$}
\end{picture}
\caption{\smaller\sf A grid property involving $P$-tiles:
assuming that the left rectangle is a $P$-tile and the bottom
sequence $(\xs_1, \xs_2)$ is normal, the top sequence $(\xs'_1,
\xs'_2)$ is normal as well.}
\label{F:Product1}
\end{figure}

\begin{lemm}
\label{L:Product1}
(Figure~\ref{F:Product1}) 
Assume that $\xs_1, \xs_2, \xs'_1, \xs'_2, \xt_0,
\xt_1, \xt_2$ be simple braids satisfying $\xt_0\xs_1=\xs'_1\xt_1$, 
$\xt_1\xs_2=\xs'_2\xt_2$, and such that $(\xs_1,\xs_2)$ and
$(\xs'_1, \xt_1)$ are normal. Then $(\xs'_1, \xs'_2)$ is normal as
well.
\end{lemm}

\begin{proof}
$(ii)$ With the notation of Figure~\ref{F:Product1}, we compute:
\begin{align*}
\HH(\xs'_1\xs'_2)
&\divel \HH(\xs'_1\xs'_2\xt_2)
&&\text{by~\eqref{E:Head3}}\\
&= \HH(\xt_0\xs_1\xs_2)
&&\text{by commutativity}\\
&= \HH(\xt_0\HH(\xs_1\xs_2))
&&\text{by~\eqref{E:Head4}}\\
&= \HH(\xt_0\xs_1)
&&\text{by the hypothesis that $(\xs_1, \xs_2)$ is normal}\\ 
&= \HH(\xs'_1\xt_1)
&&\text{by commutativity}\\
&= \xs'_1
&&\text{by the hypothesis that $(\xs'_1, \xt_1)$ is normal},
\end{align*}
so $(\xs'_1, \xs'_2)$ is normal.
\end{proof}

To go further, we need to use the duality with respect to~$\DD\nn$.
For each simple $n$-braid~$\xs$, there exists a unique simple
$n$-braid~$\xs^*$ satisfying $\xs \xs^* = \DD\nn$. Then, 
we have $\xs\DD\nn = \xs\xs^*\xs^{**} =
\DD\nn\xs^{**}$, hence $\xs^{**}$ is the conjugate~$\DD\nn\inv
\xs \DD\nn$, which is easily seen to be the image of~$\xs$ under
the flip automorphism~$\flip_\nn$ that exchanges~$\ss\ii$
and~$\ss{\nn-\ii}$ for $1 \le \ii<\nn$. 

\begin{lemm}
\label{L:Dual}
Assume that $\xs,\xt$ are simple braids. Then $(\xs,
\xt)$ is normal if and only if $\gcdl(\xs^*, \xt)=1$ holds.
\end{lemm}

\begin{proof}
We have 
$$\HH(\xs \xt) 
= \gcdl(\xs \xt, \DD\nn)
= \gcdl(\xs \xt, \xs \xs^*)
= \xs \cdot \gcdl(\xt, \xs^*),$$
hence $\HH(\xs \xt) = \xs$ is equivalent to $\gcdl(\xt, \xs^*)=1$.
\end{proof}

\begin{figure}[htb]
\begin{picture}(27,33)
\put(0.2,0.5){\includegraphics{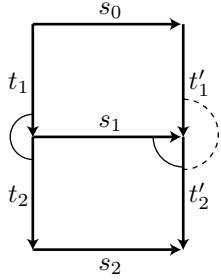}}
\put(0,23){$\xt_1$}
\put(0,8){$\xt_2$}
\put(24,23){$\xt'_1$}
\put(24,8){$\xt'_2$}
\put(12,33){$\xs_0$}
\put(12,18){$\xs_1$}
\put(12,-1){$\xs_2$}
\end{picture}
\caption{\smaller\sf Another grid property involving $P$-tiles:
assuming that the bottom rectangle is a $P$-tile and the left
sequence $(\xt_1, \xt_2)$ is normal, the right sequence $(\xt'_1,
\xt'_2)$ is normal as well.}
\label{F:Product2}
\end{figure}

\begin{lemm}
\label{L:Product2}
(Figure~\ref{F:Product2}) 
Assume that $\xs_0, \xs_1, \xs_2, \xt_1, \xt_2, \xt'_1, \xt'_2$ are
simple braids satisfying $\xt_1\xs_1=\xs_0\xt'_1$, 
$\xt_2\xs_2=\xs_1\xt'_2$, and such that $(\xt_1,\xt_2)$ and
$(\xs_1, \xt'_2)$ are normal. Then $(\xt'_1, \xt'_2)$ is normal as
well.
\end{lemm}

\begin{figure}[htb]
\begin{picture}(47,33)
\put(0.2,0.5){\includegraphics{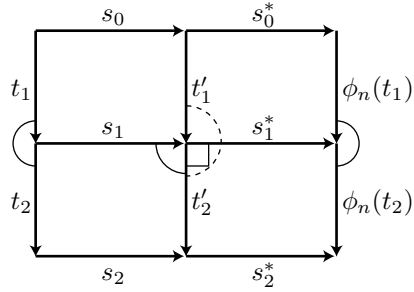}}
\put(0,23){$\xt_1$}
\put(0,8){$\xt_2$}
\put(24,23){$\xt'_1$}
\put(24,8){$\xt'_2$}
\put(12,33){$\xs_0$}
\put(12,18){$\xs_1$}
\put(12,-1.5){$\xs_2$}
\put(32,33){$\xs_0^*$}
\put(32,18){$\xs_1^*$}
\put(32,-1.5){$\xs_2^*$}
\put(44,23){$\flip_\nn(\xt_1)$}
\put(44,8){$\flip_\nn(\xt_2)$}
\end{picture}
\caption{\smaller\sf Proof of Lemma~\ref{L:Product2}: one
introduces the dual of the horizontal arrows; the
hypothesis that the left column~$(\xt_1, \xt_2)$ is normal implies
that the right column $(\flip_\nn(\xt_1),\flip_\nn(\xt_2))$ is 
normal as well, and, then, we apply
Lemma~\ref{L:Complement1} to come back to $(\xt'_1, \xt'_2)$.}
\label{F:Product2bis}
\end{figure}

\begin{proof}
Introduce $\xs_0^*$, $\xs_1^*$ and $\xs_2^*$. Then the diagram of
Figure~\ref{F:Product2bis} is commutative. As $\flip_\nn$ is an
automorphism of~$\BB\nn$---and of the group~$B_\nn$ too---the
hypothesis that $(\xt_1, \xt_2)$ is normal implies that
$(\flip_\nn(\xt_1), \flip_\nn(\xt_2))$ is normal as well. 

Now, by Lemma~\ref{L:Dual}, the hypothesis that $(\xs_1, \xt'_2)$
is normal gives $\gcdl(\xs_1^*,\xt'_2)=\nobreak1$. By
Lemma~\ref{L:Complement1}, the latter condition together with the
normality  of~$(\flip_\nn(\xt_1),
\flip_\nn(\xt_2))$ implies that $(\xt'_1,
\xt'_2)$ is normal.
\end{proof}

We deduce

\begin{figure}[htb]
\begin{picture}(72,47)
\put(0,-0.5){\includegraphics{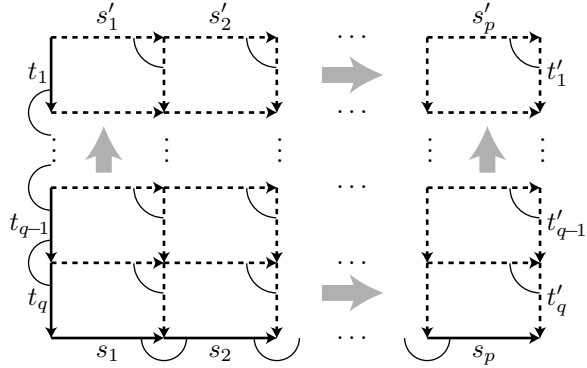}}
\put(0,37){$\xt_1$}
\put(3,17){\llap{$\xt_{\qq\!-\!1}$}}
\put(0,7){$\xt_\qq$}
\put(9,0){$\xs_1$}
\put(24,0){$\xs_2$}
\put(59,0){$\xs_\pp$}
\put(3,26){$\vdots$}
\put(18,26){$\vdots$}
\put(33,26){$\vdots$}
\put(53,26){$\vdots$}
\put(68,26){$\vdots$}
\put(41,2.5){$\ldots$}
\put(41,12.5){$\ldots$}
\put(41,22.5){$\ldots$}
\put(41,32.5){$\ldots$}
\put(41,42.5){$\ldots$}
\put(69,37){$\xt'_1$}
\put(69,17){$\xt'_{\qq\!-\!1}$}
\put(69,7){$\xt'_{\qq}$}
\put(9,44.5){$\xs'_1$}
\put(24,44.5){$\xs'_2$}
\put(59,44.5){$\xs'_\pp$}
\end{picture}
\caption{\smaller\sf Computation of the normal form of a product
by means of a grid: we start from the left and the bottom sides, and
fill the diagram using $P$-tiles, from left to right and from bottom to
top; then the top row and the right column are normal sequences as
well, and so is the sequence obtained by concatenating them.}
\label{F:Product}
\end{figure}

\begin{prop}
\label{P:Product}
$(i)$ Assume that $(\xs_1, ..., \xs_\pp)$ and $(\xt_1, ..., \xt_\qq)$
are normal sequences of simple braids. Let $\mathcal{D}$ be the grid
diagram obtained by starting from the left column $(\xt_1, ...,
\xt_\qq)$ and the bottom row $(\xs_1, ..., \xs_\pp)$ and filling the
diagram with $P$-tiles from left to right, and from bottom to top,
as shown in Figure~\ref{F:Product}. Then
every path in~$\mathcal{D}$ consisting of horizontal arrows
followed by vertical arrows corresponds to a normal sequence. 

$(ii)$ Let $\xx = \xs_1...\xs_\pp$ and $\yy=\xt_1...\xt_\qq$. Let
$(\xs'_1, ..., \xs'_\pp)$ be the top row of~$\mathcal{D}$ and
$(\xt'_1, ..., \xt'_\qq)$ is its right column. Then $(\xs'_1, ...,
\xs'_\pp, \xt'_1, ..., \xt'_\qq)$ is the normal form of~$\yy\xx$.
\end{prop}

\begin{proof}
By hypothesis, the bottom row is normal, hence, by
Lemma~\ref{L:Product1}, we inductively deduce that the $\kk$th
row from the bottom is normal. Similarly, by hypothesis, the left
column is normal, hence, by Lemma~\ref{L:Product2}, we inductively
deduce that $\kk$th column from the left is normal. Finally, the
sequence $(\xs'_1, ...,
\xs'_\pp, \xt'_1, ..., \xt'_\qq)$ is normal, as it is the concatenation
of two normal sequences and, moreover, $(\xs'_\pp, \xt'_1)$ is
normal by construction. As, by construction, the
diagram~$\mathcal{D}$ is commutative, the product of the latter
sequence is also the product of the left column and the bottom row,
\ie, it is the braid~$\yy\xx$.
\end{proof}

As in the case of the quotient, we deduce in particular rules for
computing the normal form of $\xu\xx$ and~$\xx\xu$ from that
of~$\xx$ when $\xu$ is a simple braid.

\begin{figure}[htb]
\begin{picture}(101,45)
\put(0,.3){\includegraphics{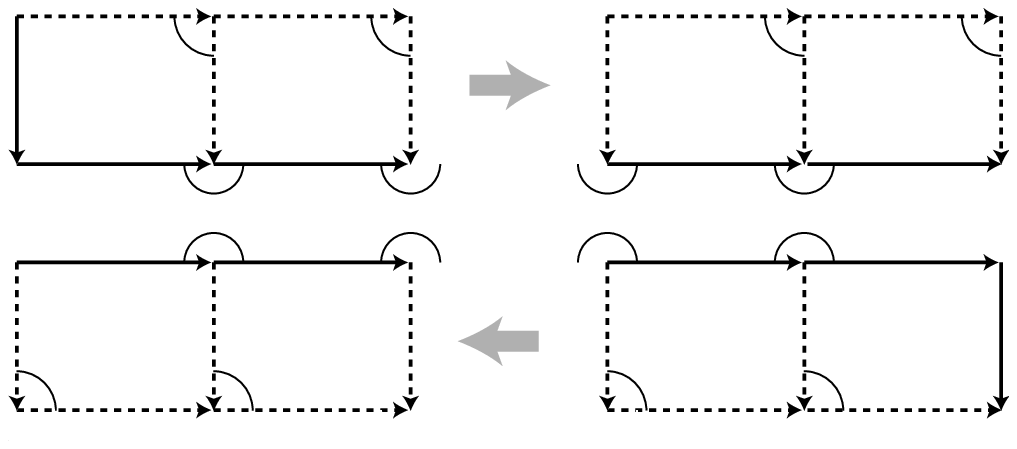}}
\put(9,0){$\xs''_1$}
\put(29,0){$\xs''_2$}
\put(68,0){$\xs''_{\pp-1}$}
\put(89,0){$\xs''_\pp$}
\put(9,20){$\xs_1$}
\put(29,20){$\xs_2$}
\put(68,20){$\xs_{\pp-1}$}
\put(89,20){$\xs_\pp$}
\put(9,26){$\xs_1$}
\put(29,26){$\xs_2$}
\put(68,26){$\xs_{\pp-1}$}
\put(89,26){$\xs_\pp$}
\put(9,45.5){$\xs'_1$}
\put(29,45.5){$\xs'_2$}
\put(68,45.5){$\xs'_{\pp-1}$}
\put(89,45.5){$\xs'_\pp$}
\put(49,3.5){$\ldots$}
\put(50,18.5){$\ldots$}
\put(49,28.5){$\ldots$}
\put(50,43.5){$\ldots$}
\put(-3.5,11){$\xu''_0$}
\put(22,11){$\xu''_1$}
\put(82,11){$\xu''_{\pp-1}$}
\put(96.5,11){$\xu''_{\pp}$}
\put(102,11){$\xu$}
\put(-2,36){$\xu$}
\put(2,36){$\xu'_0$}
\put(22,36){$\xu'_1$}
\put(82,36){$\xu'_{\pp-1}$}
\put(102,36){$\xu'_\pp$}
\end{picture}
\caption{\smaller\sf Construction of the normal
from of~$\xu\xx$ (above) and $\xx\xu$ (below) from the normal
form $(\xs_1, ..., \xs_\pp)$ of~$\xx$: start with $\xu'_0:=\xu$ (\resp 
with $\xu''_\pp:= \xu$)  and fill the diagram with $P$-tiles from left to right
(\resp from right to left); then $(\xs'_1, ..., \xs'_\pp,
\xu'_\pp)$ (\resp $(\xu''_0, \xs''_1, ..., \xs''_\pp)$) is the expected
normal form---up to removing $\xu'_\pp$ or~$\xs''_\pp$ if the
latter happen to be~$1$.}
\label{F:ProductOne}
\end{figure}

\begin{coro}
\label{C:Product}
Assume that $(\xs_1, ..., \xs_\pp)$ is the normal form
of the positive braid~$\xx$, and that $\xu$ is a simple braid. Then
the normal forms of~$\xu\xx$ and~$\xx\xu$ are the sequences
$(\xs'_1, ..., \xs'_\pp, \xu'_\pp)$ and $(\xu''_0, \xs''_1, ..., \xs''_\pp)$ 
determined by the grid diagrams of Figure~\ref{F:ProductOne}.
\end{coro}

\subsection{Application to computing normal forms}

We can now easily provide algorithms that compute the greedy
normal form and the symmetric normal form of a braid starting from
an arbitrary word that represents it. Clearly, the point is, starting 
from the greedy (\resp symmetric) normal form of a braid~$\zz$, to
be able to determine the greedy (\resp symmetric) normal form
of~$\zz\xu^{\pm1}$ when $\xu$ is a simple braid, hence in
particular a generator~$\ss\ii$. The solutions come as easy
applications of the results of the previous sections. 

\subsubsection{Computing the greedy normal form}

We recall that, for~$\xs$ a simple $\nn$-braid, $\xs^*$ denotes the
unique (simple) braid satisfying $\xs \xs^* = \DD\nn$.
Symmetrically, we denote by $\dual\xs$ the unique (simple) braid
satisfying $\dual\xs \xs=\DD\nn$. Note that, by construction, 
$(\dual\xs)^*=\xs$ holds for each simple braid~$\xs$.

\begin{figure}[htb]
\begin{picture}(120,50)
\put(0,0){\includegraphics{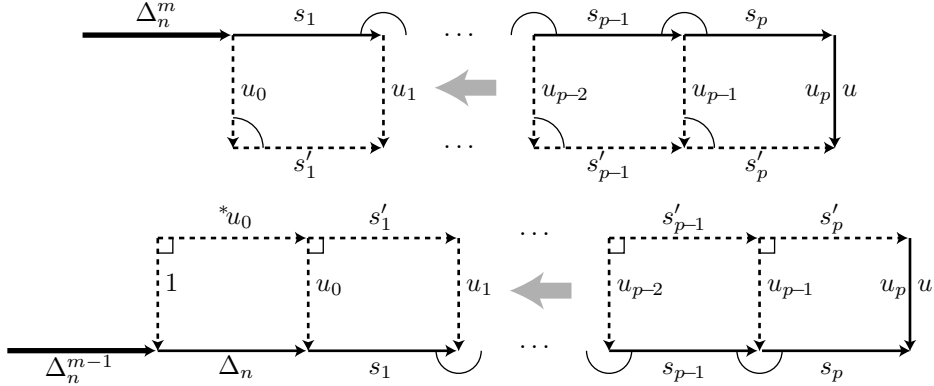}}
\put(5,-0.3){$\DD\nn^{\mm-1}$}
\put(28,0){$\DD\nn$}
\put(48,0){$\xs_1$}
\put(87,0){$\xs_{\pp\!-\!1}$}
\put(108,0){$\xs_\pp$}
\put(48,20){$\xs'_1$}
\put(87,20){$\xs'_{\pp\!-\!1}$}
\put(108,20){$\xs'_\pp$}
\put(68,3.5){$\ldots$}
\put(68,18.5){$\ldots$}
\put(28,20){$\dual{\xu_0}$}
\put(21,11){$1$}
\put(41,11){$\xu_0$}
\put(61,11){$\xu_1$}
\put(81,11){$\xu_{\pp\!-\!2}$}
\put(101,11){$\xu_{\pp\!-\!1}$}
\put(116,11){$\xu_\pp$}
\put(121,11){$\xu$}
\put(17,47){$\DD\nn^\mm$}
\put(38,47){$\xs_1$}
\put(77,47){$\xs_{\pp\!-\!1}$}
\put(98,47){$\xs_\pp$}
\put(38,27){$\xs'_1$}
\put(77,27){$\xs'_{\pp\!-\!1}$}
\put(98,27){$\xs'_\pp$}
\put(58,30){$\ldots$}
\put(58,45){$\ldots$}
\put(31,37){$\xu_0$}
\put(51,37){$\xu_1$}
\put(71,37){$\xu_{\pp\!-\!2}$}
\put(91,37){$\xu_{\pp\!-\!1}$}
\put(106,37){$\xu_\pp$}
\put(111,37){$\xu$}
\end{picture}
\caption{\smaller\sf Greedy normal form of~$\zz\xu$
(top) and~$\zz\xu\inv$ (bottom) from the greedy normal
form $(\mm; \xs_1, ..., \xs_\pp)$ of~$\zz$: starting from
$\xu_\pp:= \xu$, fill the diagram using $P$-tiles (\resp $C$-tiles)
from right to left, and adapt at the the left end to guarantee the
connection with the $\DD\nn^\mm$-factor, namely
include $\xu_0$ in the latter if needed (top), or factorize
$\DD\nn\xu_0\inv$ into~$\dual{\xu_0}$ (bottom).}
\label{F:AlgoGreedy}
\end{figure}

\begin{prop}
\label{P:AlgoGreedy}
Assume that the $\nn$-greedy normal form of~$\zz$
is $(\mm; \xs_1, ..., \xs_\pp)$ and $\xu$ is a simple $\nn$-braid. 

$(i)$ Let $\xs'_\pp$, $\xu_{\pp-1}$, ..., $\xs'_1$, $\xu_0$ be determined 
by filling the top diagram of Figure~\ref{F:AlgoGreedy}. Then the greedy 
normal form of~$\zz\xu$ is $(\mm+1; \xs'_1, ..., \xs'_\pp)$ if
$\xu_0=\DD\nn$ holds, and $(\mm; \xu_0, \xs'_1, ..., \xs'_\pp)$ otherwise.

$(ii)$ Let $\xs'_\pp$, $\xu_{\pp-1}$, ..., $\xs'_1$, $\xu_0$ be determined by
filling the bottom diagram of Figure~\ref{F:AlgoGreedy}. Then the greedy 
normal form of~$\zz\xu\inv$ is $(\mm;  \xs'_1, ..., \xs'_\pp)$ if $\xu_0=1$ holds,
and $(\mm-1; \dual{u_0}, \xs'_1, ..., \xs'_\pp)$ otherwise.
\end{prop}

\begin{proof}
$(i)$ By commutativity of the diagram, we have $\zz\xu =
\DD\nn^\mm\xu_0\xs'_1...\xs'_\pp$, so the point is to check that
the sequence $(\mm; \xu_0, \xs'_1, ..., \xs'_\pp)$, or $(\mm+1;
\xs'_1, ..., \xs'_\pp)$, is normal. Now Corollary~\ref{C:Product}
guarantees that $(\xu_0, \xs'_1, ..., \xs'_\pp)$ is normal. According
to whether $\xu_0$ equals~$\DD\nn$ or not, we integrate the
factor~$\xu_0$ in~$\DD\nn^\mm$, and we obtain a greedy normal
form, hence the greedy normal form of~$\zz\xu$ by uniqueness.

$(ii)$ By commutativity, we have $\zz\xu\inv = \DD\nn ^\mm \xu_0\inv
\xs'_1 ... \xs'_\pp = \DD\nn^{\mm-1} \dual{\xu_0} \xs'_1 ...
\xs'_\pp$, and the point is to check that the expected sequences are
greedy normal forms. Corollary~\ref{C:Complement} guarantees
that $(\xs'_1, ..., \xs'_\pp)$ is normal. So, if $\xu_0=1$
holds, $(\mm; \xs'_1, ..., \xs'_\pp)$ is a greedy normal form, hence it
is the greedy normal form of~$\zz\xu\inv$; otherwise, we observe
that $(\mm\nobreak-\nobreak1; \dual{\xu_0}, \xs'_1, ..., \xs'_\pp)$
is a greedy normal form, as $\xu_0\not=1$ implies
$\dual{\xu_0}\not=\DD\nn$, and, by Lemma~\ref{L:Dual}, 
the hypothesis $\gcdl(\xu_0, \xs'_1)=1$ is equivalent to
$(\dual{\xu_0}, \xs'_1)$ being normal since $\xu_0 =
(\dual{\xu_0})^*$ holds.
\end{proof}

\begin{exam}
\label{X:AlgoGreedy}
Let $\ww_0$ be the braid word $\mathtt{aBabacABABAbbCB}$---the
randomly chosen $4$-braid word illustrated in
Figure~\ref{F:Braid}, which will be repeatedly considered in
the sequel. Starting with $(0; \emptyset)$, which is the greedy
normal form of~$1$, and applying the algorithm of 
Proposition~\ref{P:AlgoGreedy} to the successive letters of~$\ww$,
we obtain the $4$-greedy normal form of the prefixes
of~$\ww$, namely (in braid words form):

$0: \e \gives (0; \emptyset)$

$1: \mathtt{a} \gives (0; \mathtt{a})$

$2: \mathtt{aB} \gives (-1; \mathtt{abcb,ba})$

$3: \mathtt{aBa} \gives (-1; \mathtt{abcb,ba,a})$

$4: \mathtt{aBab} \gives (-1; \mathtt{abcb,ba,ab})$

$5: \mathtt{aBaba} \gives (0; \mathtt{a,ab})$

$6: \mathtt{aBabac} \gives (0; \mathtt{a,abc})$

\dots

$14:  \mathtt{aBabacABABAbbC} \gives (-2;
\mathtt{ac, abcb, bcba, ab})$

$15:  \mathtt{aBabacABABAbbCB} \gives (-2;
\mathtt{ac, abcb, bcba, a})$.

\noindent Thus the greedy normal form of the braid represented
by~$\ww$ is the sequence

$(-2; \mathtt{ac, abcb, bcba, a})$, 

\noindent\ie,
in permutation form, 

$(-2; (2,1,4,3), (2,4,3,1), (4,1,3,2), (2,1,3,4))$

\noindent---this is the greedy normal
form of Example~\ref{X:NF}. As the initial word~$\ww_0$ contains
$15$~letters, computing its greedy normal form entails
$15$~applications of Proposition~\ref{P:AlgoGreedy}. However, one
can speed up the process by gathering adjacent letters that
together represent a simple braid: for instance, $\mathtt{abac}$
represents a simple braid, so steps~$3$ to~$6$ above
can be gathered into a single step corresponding to multiplying
by the simple braid~$\mathtt{abac}$. By Corollary~\ref{C:Greedy},
we deduce that the braid word~$\ww_0$ does not represent~$1$
in~$B_4$.
\end{exam}

As for complexity analysis, two cases are to be considered. If the
braid index~$\nn$ is fixed, and, for practical implementations, has a
small value, say $\nn\le6$, then one can precompute the
tables of the binary operations $(\xs,\xt)\mapsto \HH(\xs\xt)$,
$(\xs,\xt)\mapsto \HH(\xs\xt)\inv\xs\xt$, and $(\xs,\xt)\mapsto
\gcdl(\xs,\xt)$, in which case the determination of the greedy
normal form for a braid word of length~$\ell$ has
complexity~$O(\ell^2)$, and is easy and quick in practice. Otherwise,
there are too many simple $\nn$-braids to store all results, and one
has to compute the values of~$\HH(\xs\xt)$ and~$\gcdl(\xs,\xt)$
locally. As explained in~\cite[Chapter 9]{Eps},
this is essentially a sorting process, and, therefore, each such
computation can be done in time~$O(\nn\log\nn)$, resulting in a
global time complexity~$O(\ell^2\nn\log\nn)$ for the computation
of the normal form for an $\nn$-braid word of length~$\ell$.

\subsubsection{Computing the symmetric normal form}

We now consider the symmetric normal form of
Section~\ref{S:Symmetric}, and show how to compute the symmetric
normal form of~$\zz\xu$ and $\zz\xu\inv$ form that of~$\zz$
by using convenient grid diagrams. The method is similar to that for
the greedy normal, but a bit more care is needed for the transition
between the numerator and the denominator in the case of a
product.

\begin{figure}[htb]
\begin{picture}(140,51)
\put(-10,-0.5){\includegraphics{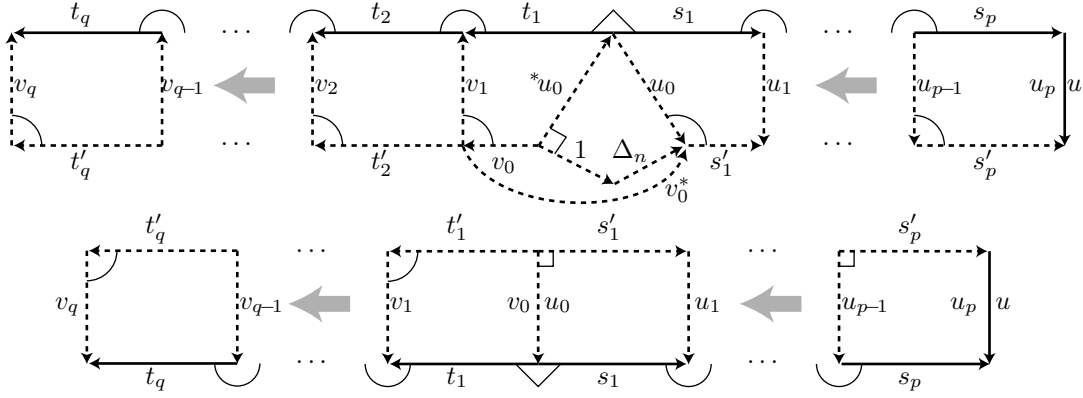}}
\put(9,0){$\xt_\qq$}
\put(49,0){$\xt_1$}
\put(69,0){$\xs_1$}
\put(109,0){$\xs_\pp$}
\put(9,20){$\xt'_\qq$}
\put(49,20){$\xt'_1$}
\put(69,20){$\xs'_1$}
\put(109,20){$\xs'_\pp$}
\put(29,2.7){$\ldots$}
\put(29,17.5){$\ldots$}
\put(89,2.7){$\ldots$}
\put(89,17.5){$\ldots$}
\put(-3,10){$\xv_\qq$}
\put(21.5,10){$\xv_{\qq\!-\!1}$}
\put(41.5,10){$\xv_1$}
\put(57,10){$\xv_0$}
\put(62,10){$\xu_0$}
\put(82,10){$\xu_1$}
\put(101.5,10){$\xu_{\pp\!-\!1}$}
\put(116,10){$\xu_\pp$}
\put(122,10){$\xu$}

\put(-1,48.5){$\xt_\qq$}
\put(39,48.5){$\xt_2$}
\put(59,48.5){$\xt_1$}
\put(79,48.5){$\xs_1$}
\put(119,48.5){$\xs_\pp$}
\put(-1,28.5){$\xt'_\qq$}
\put(39,28.5){$\xt'_2$}
\put(78,25){$\xv_0^*$}
\put(84,28.5){$\xs'_1$}
\put(119,28.5){$\xs'_\pp$}
\put(19,46.7){$\ldots$}
\put(99,46.7){$\ldots$}
\put(19,31.7){$\ldots$}
\put(99,31.7){$\ldots$}
\put(-8.5,39){$\xv_\qq$}
\put(11,39){$\xv_{\qq\!-\!1}$}
\put(31.5,39){$\xv_2$}
\put(51.5,39){$\xv_1$}
\put(76,39){$\xu_0$}
\put(55,28.5){$\xv_0$}
\put(66,30){$1$}
\put(71,30){$\DD\nn$}
\put(60,39){$\dual{\xu_0}$}
\put(91.5,39){$\xu_1$}
\put(111.5,39){$\xu_{\pp\!-\!1}$}
\put(126.5,39){$\xu_\pp$}
\put(131.5,39){$\xu$}
\end{picture}
\caption{\smaller\sf Symmetric normal form of~$\zz\xu$
(top) and~$\zz\xu\inv$ (bottom) from the symmetric
normal form $(\xt_1, ..., \xt_\qq; \xs_1, ..., \xs_\pp)$
of~$\zz$: start from $\xu_\pp := \xu$, and fill the diagram 
from right to left using $P$-tiles  (\resp
$C$-tiles, then $P$-tiles); in the case of the product, the transition 
is as follows: having got~$\xu_0$, we find $\dual{\xu_0}$ using a
$C$-tile, then let $(\xv_0, \xv_1)$ be the normal form of $\dual{\xu_0}\xt_1$
using a $P$-tile, and
continue with $\xv_1, \xt_2$, \etc.; in the case of the
quotient, we simply put $\xv_0:=\xu_0$, and continue
with $\xv_0, \xt_1$, \etc}
\label{F:AlgoSymmetric}
\end{figure}

\begin{prop}
\label{P:AlgoSymmetric}
Assume that the symmetric normal form of~$\zz$
is the double sequence $(\xt_1, ..., \xt_\qq; \xs_1, ..., \xs_\pp)$, and
$\xu$ is a simple braid. 

$(i)$ Let $\xs'_\pp$, $\xu_{\pp-1}$, ..., $\xs'_1$, $\xu_0$, $\xv_0$,
$\xv_1$, $\xt'_2$, $\xv_2$,  ... , $\xt'_\qq$, $\xv_\qq$ be
determined by filling the top diagram of
Figure~\ref{F:AlgoSymmetric}. Then the symmetric normal 
form of~$\zz\xu$ is $(\xt'_2, ..., \xt'_\qq, \xv_\qq; \xv_0^*,  \xs'_1, ..., \xs'_\pp)$, 

$(ii)$ Let $\xs'_\pp$, $\xu_{\pp-1}$, ..., $\xs'_1$, $\xu_0=\xv_0$,
$\xv_1$, $\xt'_2$, ... , $\xt'_\qq$, $\xv_\qq$ are determined by
filling the bottom diagram of Figure~\ref{F:AlgoSymmetric}. Then the 
symmetric normal form of~$\zz\xu\inv$ is $(\xt'_1, ...,
\xt'_\qq, \xv_\qq;  \xs'_1, ..., \xs'_\pp)$ , 
\end{prop}

\begin{proof}
$(i)$ By commutativity of the diagram, we have $\zz\xu =
\xv_\qq\inv \xt'_\qq{}\inv ... \xt'_2{}\inv \xv_0^* \xs'_1 ... \xs'_\pp$, and
the point is to show that the double sequence $(\xt'_2, ..., \xt'_\qq, \xv_\qq; 
\xv_0^*, \xs'_1, .., \xs'_\pp)$ is a symmetric normal form. 
Corollary~\ref{C:Product} guarantees that the sequences $(\xu_0, \xs'_1, ..., \xs'_\pp)$
and $(\xv_0, \xt'_2, ..., \xt'_\qq, \xv_\qq)$ are normal. There
remain two points to check, namely that $\gcdl(\xt'_2,
\xv_0^*)$ is~$1$, and that $(\xv_0^*,\xs'_1)$ is normal.
For the first relation, Lemma~\ref{L:Product1} implies that $(\xv_0,
\xt'_2)$ is normal, which is equivalent to
$\gcdl(\xt'_2,
\xv_0^*)=\nobreak1$ by Lemma~\ref{L:Dual}.  

As for the second relation, \ie, for the normality of~$(\xv_0^*,\xs'_1)$,
 by Lemma~\ref{L:Dual} again, it is equivalent to
$\gcdl(\xv_0^{**},\xs'_1)=1$, hence to
$\gcdl(\xv_0,\flip_\nn\inv(\xs'_1))=1$. Now  we have
$\xv_0\xv_1= \dual{\xu_0}\xt_1$, and $\xs_1\xu_1 =
\xu_0\xs'_1$, hence 
$$\dual{\xu_0} \xs_1 \xu_1 
= \DD\nn \xs'_1
= \flip_\nn\inv(\xs'_1) \DD\nn
= \flip_\nn\inv(\xs'_1) \dual{\xu_1} \xu_1,$$
hence $\flip_\nn\inv(\xs'_1) \dual{\xu_1} = \dual{\xu_0}\xs_1$.
Assume $\xs$ is a simple left divisor of~$\xv_0$
and~$\flip_\nn\inv(\xs'_1)$. Then, {\it a fortiori}, we have
$\xs\divel \xv_0 \xv_1$ and $\xs \divel
\flip_\nn\inv(\xs'_1)\dual{\xu_1}$, hence, by the above
computations, $\xs \divel \dual{\xu_0} \xt_1$ and $\xs \divel
\dual{\xu_0} \xs_1$, and, therefore, $\xs \divel \dual{\xu_0}$ as,
by hypothesis, $\gcdl(\xs_1, \xt_1)= 1$ holds. Now, by hypothesis,
$(\xu_0, \xs'_1)$ is normal, hence, by Lemma~\ref{L:Dual}, we have
$\gcdl(\xu_0^*, \xs'_1)=1$, and therefore
$\gcdl(\flip_\nn\inv(\xu_0^*), \flip_\nn\inv(\xs'_1)) =1$, \ie,
$\gcdl(\dual{\xu_0}, \flip_\nn\inv(\xs'_1))=1$. So we must have
$\xs=1$, implying $\gcdl(\xv_0,\flip_\nn\inv(\xs'_1))=1$, which
was seen above to be equivalent to $(\xv_0^*,\xs'_1)$ being normal.
So the proof is complete.

$(ii)$ The argument for the quotient is similar. By commutativity of
the diagram, we have $\zz\xu\inv = \xv_\qq\inv
\xt'_\qq{}\inv ... \xt'_1{}\inv \xs'_1... \xs'_\pp$, and, once
again, the point is to check that the double sequence
$(\xt'_1, ..., \xt'_\qq, \xv_\qq; \xs'_1, ..., \xs'_\pp)$ is a
symmetric normal form. Corollary~\ref{C:Complement} 
implies that $(\xs'_1, ..., \xs'_\pp)$ is normal, and
Corollary~\ref{C:Product} implies that $(\xt'_1, ..., \xt'_\qq,
\xv_\qq)$ is normal, so the only remaining point to check is
$\gcdl(\xs'_1, \xt'_1)=1$. Now assume that $\xs$ is a simple left
divisor of~$\xs'_1$ and~$\xt'_1$. Then, {\it a fortiori}, we have
$\xs \divel \xs'_1\xu_1 = \xu_0\xs_1$ and $\xs \divel \xt'_1 \xv_1
= \xu_0 \xt_1$. As $\gcdl(\xs_1, \xt_1)=1$ holds by hypothesis, we
deduce $\xs \divel \xu_0$, and, finally, $\xs=1$ as, by construction,
we have $\gcdl(\xu_0, \xs'_1)=1$. Hence $\gcdl(\xs'_1, \xt'_1)=1$
holds.
\end{proof}

\begin{exam}
\label{X:AlgoSymmetric}
Let $\ww_0$ be the braid word $\mathtt{aBabacABABAbbCB}$ again.
Applying the algorithm of  Proposition~\ref{P:AlgoSymmetric} to the
successive letters of~$\ww$ leads to the symmetric normal forms
(here in braid word form):

$0: \e \gives (\emptyset; \emptyset)$

$1: \mathtt{a} \gives (\emptyset; \mathtt{a})$

$2: \mathtt{aB} \gives (\mathtt{ab}; \mathtt{ba})$

$3: \mathtt{aBa} \gives (\mathtt{ab}; \mathtt{ba,a})$

$4: \mathtt{aBab} \gives (\mathtt{ab}; \mathtt{ba,ab})$

$5: \mathtt{aBaba} \gives (\emptyset; \mathtt{a,ab})$

$6: \mathtt{aBabac} \gives (\emptyset; \mathtt{a,abc})$

\dots

$14:  \mathtt{aBabacABABAbbC} \gives (\mathtt{ab, bacb};
\mathtt{bcba, ab})$

$15:  \mathtt{aBabacABABAbbCB} \gives (\mathtt{ab, bacb};
\mathtt{bcba, a})$.

\noindent Thus the symmetric normal form of the braid represented
by~$\ww_0$ is the double sequence  $(\mathtt{ab, bacb};
\mathtt{bcba, a})$, \ie, in permutation form, 

$((2,3,1,4), (3,4,1,2); (4,1,3,2), (2,1,3,4))$

\noindent---this is the symmetric
normal form of Example~\ref{X:Symmetric}. As in the case of the
greedy normal form, we observe that the steps corresponding to a
single simple factor can be gathered. We deduce from
Corollary~\ref{C:Symmetric} that $\ww_0$ does not represent~$1$
in~$B_4$.
\end{exam}

The complexity analysis is the same as in the case of the greedy
normal form: for a fixed braid index~$\nn$, the algorithm of
Proposition~\ref{P:AlgoSymmetric} is quadratic in the length of the
initial word; when $\nn$ is not fixed, a multiplicative
factor~$\nn\log\nn$ has to be inserted.

\begin{rema}
The diagrams of Figures~\ref{F:AlgoGreedy}
and~\ref{F:AlgoSymmetric} make it clear that the greedy and
symmetric normal forms satisfy the fellow traveler property
of~\cite{Eps} and therefore are connected with an automatic
structure on the braid group~$B_\nn$.

Also, let us point here that the greedy and symmetric normal 
forms exist for a class of structures that is much wider than braid
groups, namely the so-called Garside groups of~\cite{Dfx, Dgk}, and
even more: thin groups of fractions of~\cite{Dgo}, Garside
categories of~\cite{Kra, DiM, Bes}---all eligible for the grid properties
described above.

Finally, we mention the existence of alternative normal
forms~\cite{Bur, Dho} whose computation has the same complexity as
the greedy and symmetric normal forms, but which rely on different
bases and are connected with the braid order alluded to in
Section~\ref{S:Handle} below. So far, there seems to be no reason to
expect these normal forms to be more suitable for practical
applications than those described above.
\end{rema}

\section{Direct solutions}
\label{S:Direct}

Using a normal form is not the only way for solving the braid isotopy
problem. Besides the solutions of Section~\ref{S:Normal}, there exist
alternative solutions directly deciding whether a given braid
word~$\ww$ represents the trivial braid~$1$ or not. 

Using such solutions entails working with arbitray
braid words. As discussed in~\cite{Dgw}, this option makes the
computation of product and inverse obvious (merely concatenating
or reversing words), at the expense of making equivalence not
obvious, whereas the option of using a normal form and restricting
to normal words makes equivalence obvious (a mere equality), but
makes the algebraic operations of product and inverse less obvious,
as normalization processes such as those described in
Propositions~\ref{P:Product} and~\ref{P:Complement} are needed.

Here we describe three direct solutions to the braid isotopy problem,
namely two syntactic solutions based on some word rewrite systems,
and one geometric method due to I.\,Dynnikov which consists in
attributing integral coordinates to every braid.

\subsection{Word redressing}
\label{S:Redressing}

Word redressing (also called reversing in literature) is a simple
syntactic transformation that consists in pushing the positive  
letters~$\ss\ii$ in one direction and the negative letters~$\sss\ii$ in
the other direction, until a word of the form ``all positive, then all
negative'' is obtained. This leads to solutions of the word problem
that admit a quadratic complexity. The underlying theory is Garside
theory again.

\subsubsection{Description}

Redressing a braid word consists in looking for the subwords of the
form~$\sss\ii\ss\jj$, \ie, one negative letter followed by a positive
letter, and transforming them into equivalent patterns consisting of
positive letters followed by negative letters, \ie, replacing
negative--positive subwords with equivalent positive--negative words.

\begin{defi}
\label{D:Redr}
\cite{Dff, Dgp}
Assume that $\ww, \ww'$ are braid words. We say that $\ww$ is {\it
right redressible} to~$\ww'$ in one step if $\ww'$ is obtained
from~$\ww$ either by deleting some subword of the
form~$\sss\ii\ss\ii$, or by replacing some subword of the
form~$\sss\ii\ss\jj$ with $\vert\ii-\jj\vert\ge2$ by $\ss\jj\sss\ii$,
or by replacing some subword of the
form~$\sss\ii\ss\jj$ with $\vert\ii-\jj\vert=1$ by
$\ss\jj\ss\ii\sss\jj\sss\ii$. We say that $\ww\redrr\ww'$ holds if
there is a finite sequence $\ww_0=\nobreak\ww$, ..., $\ww_\NN=\ww'$
such that $\ww_\kk$ is right redressible to~$\ww_{\kk+1}$ in one step
for each~$\kk$.
\end{defi}

By construction, the words that are terminal for redressing are the
words that include no subword of the form~$\sss\ii\ss\jj$, \ie, the
words of the form~$\uu\vv{}\inv$ with $\uu,\vv$ positive braid
words. 

\begin{theo}
\cite{Dff}
\label{T:Redressing}
For each braid word~$\ww$, there exist two unique positive
braid words~$\uu,\vv$ such that $\ww \redrr \uu\vv{}\inv$ holds.
\end{theo}

A new solution to the braid isotopy problem follows: 

\begin{coro}
\label{C:Redressing}
A braid word~$\ww$ represents~$1$ in the braid group if and
only if, denoting by~$\uu$ and~$\vv$ the positive words for which
$\ww \redrr \uu\vv{}\inv$ holds, we have $\vv{}\inv\uu \redrr
\e$, where $\e$ denotes the empty word.
\end{coro}

\begin{exam}
\label{X:Redressing}
Let us start again with the braid
word~$\ww_0 = \mathtt{aBabacABABAbbCB}$ of
Examples~\ref{X:AlgoGreedy} and~\ref{X:AlgoSymmetric}. Owing to
Corollary~\ref{C:Redressing}, we decide whether $\ww$ is trivial
or not by redressing~$\ww$ to a word of the form~$\uu\vv{}\inv$
with $\uu, \vv$ positive words, then redressing $\vv{}\inv\uu$ again,
and looking whether we finally obtain the empty word. In the
current case, selecting at each step the leftmost pattern of the
form~$\sss\ii\ss\jj$ (underlined in the words
below), the successive words are as follows

$0: \mathtt{a\underline{Ba}bacABABAbbCB}$

$1: \mathtt{aabA\underline{Bb}acABABAbbCB}$

$2: \mathtt{aab\underline{Aa}cABABAbbCB}$

$3: \mathtt{aabcABAB\underline{Ab}bCB}$

$4: \mathtt{aabcABA\underline{Bb}aBAbCB}$

\ldots

$11: \mathtt{aabcbABBACB}$.

\noindent At this point, we switch the positive and negative
subword, and redress again:

$12: \mathtt{ABBAC\underline{Ba}abcb}$

$13: \mathtt{ABBA\underline{Ca}bABabcb}$

$14: \mathtt{ABB\underline{Aa}CbABabcb}$

\ldots

$37: \mathtt{cbaa\underline{Bb}cBCABBA}$.

$38: \mathtt{cbaacBCABBA}$.

\noindent The latter word is not empty: we conclude that
$\ww_0$ does not represent~$1$~in~$B_4$.
\end{exam}

\subsubsection{Explanation}

Clearly, redressing replaces a braid word with an equivalent braid word,
and it is easy to show that, if $\ww$ redresses to a word of the form
$\uu\vv{}\inv$ with $\uu,\vv$ positive, then $\uu$ and $\vv$ are uniquely
determined. Now the problem is that, as Example~\ref{X:Redressing} shows, 
redressing may increase the length of the words, and it is not obvious that
the process terminates. The point is that, if $\ww$ has the form $\ww_1\inv
\ww_2$ where $\ww_1, \ww_2$ are positive words representing simple
braids, then $\ww$ redresses to some word $\uu\vv{}\inv$ where $\uu$
and~$\vv$ are positive and again represent simple braids. Thus, with respect
to an enhanced alphabet containing all positive braid words representing
simple braids, the length does not increase under redressing, and this
leads to the termination result. Actually,
what redressing does is to compute the right lcm in the braid monoid, and
so the properties behing redressing are the Garside theory.

\subsubsection{Discussion}

The advantage of the redressing method is the simplicity of its
implementation: there is a unique syntactic operation, 
involving length~$2$ subwords of the considered word only.
From that point of view, the method is more easily implemented than
the greedy or symmetric normal form. At a
theoretical level, both methods are essentially equivalent: there exist
positive constants~$C, C'$ such that, for each $\nn$-braid
word~$\ww$, if $\NN_1(\ww)$ (\resp $\NN_2(\ww)$) denotes the
number of braid relations needed to put~$\ww$ into a greedy normal
form (\resp to redress~$\ww$), then we have
$C\NN_1(\ww) \le \NN_2(\ww) \le C'\NN_1(\ww)$. However,
by representing simple braids by permutations and using 
fast sorting algorithms to compute the normal
form as explained in~\cite[Chapter 9]{Eps}, one presumably obtains
a more efficient algorithm.

According to Theorem~\ref{T:Redressing}, the redressing method
starting from a braid word~$\ww$ yields a final braid word~$\ww'$
such that $\ww\equiv\e$ is equivalent to $\ww'=\e$. However, in
general, $\ww'\equiv\ww$ fails: because
of the exchange of the negative and positive factors between the two
passes, $\ww'$ is only equivalent to a conjugate of~$\ww$. Now, it is easy to describe a variant of the
redressing method that avoids the median conjugation. Indeed, let 
{\it left redressing}, denoted $\ww\redrl\ww'$, be the symmetric 
counterpart to (right)
redressing consisting in replacing each pattern~$\ss\ii\sss\jj$ with
an equivalent negative--positive word. Formally, we may define
$\ww\redrl\ww'$ to mean
$\widetilde\ww\redrr \widetilde{\ww'}$, where
$\widetilde\ww$ is the word obtained from~$\ww$ by
reversing the order of the letters.

\begin{coro}
\cite{Dff}
\label{C:RedressingBis}
A braid word~$\ww$ represents~$1$ in the braid group if and only if, denoting
by~$\uu$ and~$\vv$ the positive words for which $\ww
\redrr \uu\vv{}\inv$ holds, we have $\uu\vv{}\inv \redrl \e$.
\end{coro}

So the method is the same as in Corollary~\ref{C:Redressing},
with the only difference that we do not change the word
obtained at the end the first redressing pass, but instead continue
with left redressing. The criterion remains the same, namely that the
initial word~$\ww$ is trivial if and only if the final word~$\ww''$ is
empty. The advantage of this variant is that $\ww''$ is equivalent
to~$\ww$, and, moreover, it gives a fractionary decomposition
of~$\cl\ww$ which is geodesic, \ie, has minimal length among all
fractionary expressions of~$\cl\ww$. This implies that
the negative and the positive parts of~$\ww''$ must be equivalent to
the two components of the symmetric normal form
of~$\cl\ww$, although they need not be in normal form in general.
Assuming that an algorithm for computing the normal form
of a positive braid is available, one obtains in this way an alternative way 
for computing the symmetric normal form of an arbitrary braid that is
more simply implemented than the method of Proposition~\ref{P:AlgoSymmetric}: 
perform double redressing, then put the numerator and the denominator 
in normal form. 

\begin{exam}
Starting with~$\mathtt{aBabacABABAbbCB}$ once more, the first
$11$~steps of the redressing algorithm are as in
Example~\ref{X:Redressing}, but, then, we appeal to left redressing,
and the sequel is different. We find:

$0: \mathtt{a\underline{Ba}bacABABAbbCB}$

\ldots

$11: \mathtt{aabc\underline{bA}BBACB}$

$12: \mathtt{aab\underline{cA}BabBBACB}$

$13: \mathtt{aa\underline{bA}cBabBBACB}$

$14: \mathtt{a\underline{aA}BabcBabBBACB}$

\ldots

$39: \mathtt{BACBBAcbaac\underline{bB}}$.

$40: \mathtt{BACBBAcbaac}$.

\noindent Once again, the latter word is not empty, and we
conclude that $\ww_0$ does not represent~$1$ in~$B_4$. In
addition, we obtain that $\mathtt{BACBBAcbaac}$ is a shortest
expression of~$\cl{\ww_0}$ as a negative--positive fraction, which is
coherent with Example~\ref{X:Symmetric}, where it was shows that
the symmetric normal form of~$\cl{\ww_0}$ is $(\mathtt{ab, bacb},
\mathtt{bcba, a})$: indeed,
$(\mathtt{ab, bacb})$ is the normal form of~$\mathtt{abbcab}$,
and $(\mathtt{bcba, a})$ is the normal form of~$\mathtt{cbaac}$.
\end{exam}

As a final remark, let us observe that, because redressing is efficient
at computing lcm's and complements in the braid monoid, it also
provides an easy way to compute gcd's and, from there, normal forms.
Though not as efficient as those based on quick sorting, the algorithms
based on word redressing are more easily implemented than the 
latter, and they are convenient for small and medium size braid words.

\subsection{Handle reduction}
\label{S:Handle}

Handle reduction is another syntactic braid word transformation
which, like redressing, consists in iterating some basic word
transformation and concluding that the initial braid word
represents~$1$ if and only if the final word is empty. 
The basic transformation step is more complicated that the one
involved in redressing but the number of steps is much lower, and
the method turns out to be extremaly efficient in practice. The
underlying structure behind handle reduction is a linear
ordering of braids, which pilots the reduction process and
heuristically explains its efficiency.

\subsubsection{Description}

Handle reduction is an extension of free reduction. The latter 
consists in iteratively deleting patterns
of the form $x x\inv$ or $x\inv x$. Handle reduction involves not
only patterns of the form~$\ss i \sss i$ or $\sss i \ss i$, but also
more general patterns of the form $\ss i \dots
\sss i$ or $\sss i \dots \ss i$ with intermediate letters between the
letters~$\ss i$ and~$\sss i$.

\begin{defi}
$(i)$ A {\it $\ss i$-handle} is a braid word of the form
\begin{equation}\label{E:Handle}
w = \s_i^e \, w_0 \, \s_{i+1}^d \, w_1 \, \s_{i+1}^d
\dots
\s_{i+1}^d \, w_m \, \s_i^{-e},
\end{equation}
with $e, d = \pm 1$, $m \ge 0$, and $w_0, ..., w_m$
containing no~$\ssss j$ with $j \le i+1$. 
Then the {\it  reduct} of~$w$ is defined to be
\begin{equation}
w' = w_0 \, \s_{i+1}^{-e}\s_i^d\s_{i+1}^e \, w_1 \,
\s_{i+1}^{-e}\s_i^d\s_{i+1}^e \dots
\s_{i+1}^{-e}\s_i^d\s_{i+1}^e \, w_m,
\end{equation}
\ie, we delete the initial and final letters~$\ssss i$, and 
we replace each letter~$\ssss{i+1}$ with
$\s_{i+1}^{-e}\ssss i\s_{i+1}^e$. 

$(ii)$ We say that a braid word~$\ww$ is {\it reduced} if it
contains no $\ss\ii$-handle, where $\ss\ii$ is the generator with
minimal index occurring in~$\ww$.
\end{defi}

\begin{figure}[htb]
\setlength{\unitlength}{1mm}
\begin{picture}(125,20)
\put(0,0){\includegraphics[scale=0.9]{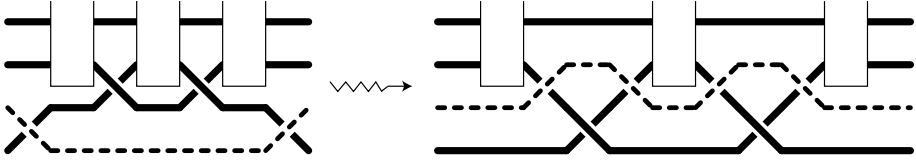}}
\end{picture}
\caption{\smaller\sf Reduction of a handle, here a $\ss1$-handle: 
the dashed strand has the shape of a handle, and reduction consists
in pushing that strand so that it skirts above the next crossings
instead of below.}
\label{F:Handle}
\end{figure}

A braid word of the form $\ss i \sss i$ or $\sss i \ss i$ is a handle,
and reducing it means deleting it, so handle reduction generalizes
free reduction.  As illustrated in Figure~\ref{F:Handle}, reducing
a handle yields an equivalent braid word. So, as in the case of
free reduction, if there is a reduction sequence from a braid
word~$w$ to the empty word, \ie, a sequence $w = w_0, 
..., w_N = \e$ such that, for each~$k$, the
word~$w_{k+1}$ is obtained from~$w_k$ by replacing
some handle of~$w_k$ by its reduct, then $w$ represents~$1$
 in the braid group. The
point is that the converse is also true.

\begin{theo}
\cite{Dfo}
\label{P:HRed}
For every braid word~$\ww$, every sequence of handle
reductions from~$\ww$ leads in finitely many steps to a reduced
word~$\ww'$. Moreover, a reduced word~$\ww'$ represents~$1$ if
and only if it is empty.
\end{theo}

We obtain a new solution to the braid isotopy problem:

\begin{coro}
\label{C:Handle}
A braid word~$\ww$ represents~$1$ in the braid group if and only
if any sequence of handle reductions starting from~$\ww$
terminates with the empty word.
\end{coro}

\begin{exam}
\label{X:Handle}
Consider $\ww_0 = \mathtt{aBabacABABAbbCB}$ again. Choosing to
reduce the leftmost handle at each step and underlying it, we
succesively obtain:

$0: \mathtt{aBab\underline{acA}BABAbbCB}$

$1: \mathtt{aBa\underline{bcB}ABAbbCB}$

$2: \mathtt{aB\underline{aCbcA}BAbbCB}$

$3: \mathtt{aBCBa\underline{bcB}AbbCB}$

$4: \mathtt{aBCB\underline{aCbcA}bbCB}$

$5: \mathtt{aBCBCBabcbbCB}$.

\noindent The latter word contains no $\ss1$-handle, so it is
reduced, and it is not empty, so we conclude that $\ww_0$ does not
represent~$1$ in~$B_4$.
\end{exam}

\subsubsection{Explanation}

Two different structures lie behind handle reduction, namely
Garside's theory that was already involved on the previous solutions,
and, in addition, some linear order that is compatible with left
multiplication on~$B_\nn$. It can be seen that each handle reduction
is essentially a composition of redressing steps, but the main
difference with the algorithms of Sections~\ref{S:Normal}
and~\ref{S:Redressing} is that, here, we do not perform all
redressing steps systematically, but only some of them
according to a general strategy provided by the underlying braid
order. This should make it natural why the handle reduction method
is, in practice, much more efficient than the redressing method and
than the greedy and symmetric normal form methods ($5$
\vs $40$ steps in our example).

\subsubsection{Discussion}

A braid word may contain many handles, so building an actual
algorithm requires to fix a strategy prescribing in which
order the handles will be reduced. In Example~\ref{X:Handle},
we chose to reduce the leftmost handle, but more
efficient strategies exist. As can be expected, the most
efficient ones use a divide-and-conquer trick. Although the only
upper bound for space and time complexity proved so far is exponential, 
handle reduction is extremely efficient in practice, as show the statistics
of~\cite{Dgw}. Also, reduction being a local
procedure, the amount of memory needed to implement it is what is
needed to just store the braid under reduction. So, using arbitrary 
words together with handle reduction instead of normal words could
be specially interesting when the computing resources are limited. 

It was mentioned above that the symmetric normal form, as well as the
redressing method, yield fractionary expressions of minimal length.
Such fractionary expressions need not be expressions of
minimal length: there may be shorter expressions that are not
fractions, \ie, in which all negative letters are not gathered in one
block and all positive letters in another block. The general problem
of finding the shortest expression of a braid is difficult: its $\Bi$
version is known to be $NP$-complete~\cite{PaR}. Although no actual
result is proved, it has been observed that handle reduction is good at
providing short expressions. Actually, the definition of handle
reduction is not symmetric and the left side plays a distinguished
r\^ole. One can compensate this lack of symmetry by defining {\it
iterated handle reduction} as follows: starting with~$\ww$, we
reduce~$\ww$ to~$\ww'$, then flip~$\ww'$ using~$\flip_\nn$,
reduce it, and flip the result. Equivalently, after one handle
reduction, one performs the symmetric operation in which the right
side is distinguished. By doing so, and possibly iterating the process,
one empirically obtains short expressions. It is conjectured by the
author, as well as by A.\,Miasnikov and A.\,Ushakov, that there
might exist a constant~$C$ such that applying the previous process
to any braid word~$\ww$ leads to a final word of length at
most~$C\ell_{\mathrm{min}}(\cl\ww)$, where, for~$\xx$ a braid, 
$\ell_{\mathrm{min}}(\xx)$ denotes the length of the shortest
word~$\xx$.

\begin{exam}
\label{X:HandleBis}
Starting once more with $\ww_0 = \mathtt{aBabacABABAbbCB}$,
iterated handle reduction from~$\ww_0$ leads in 3~steps to the
word $\mathtt{acBCCBa}$, which happens to be a geodesic
representative of the braid~$\cl{\ww_0}$.
\end{exam}

\subsection{Dynnikov coordinates}
\label{S:Dynnikov}

We conclude with still another solution to the braid isotopy
problem, namely the one provided by the so-called Dynnikov
coordinates. This solution relies on a completely different approach
steming from geometry and deep results by L.\,Mosher about
the existence of an automatic structure for all mapping class
groups~\cite{Mos}.

\subsubsection{Description}

The principle consists in associating with every $\nn$-braid word a
sequence of $2\nn$~integers, that can be thought of as
coordinates for the braid. This coordinization is faithful in that
two braid words receive the same coordinates if and only if they
represent the same braid.

\begin{defi}
$(i)$ For $\xx$ in~Ê$\Int$, write $\Pos\xx$ for~$\max(0,\xx)$,
and $\Neg\xx$ for~$\min(\xx,0)$. Let~$\FDyn,
\GDyn: \Int^4 \to \Int^4$ be defined by $\FDyn =
(\FDyn_1, ..., \FDyn_4)$, $\GDyn=(\GDyn_1, ...,
\GDyn_4)$ with
\begin{equation}
\label{E:Dynnikov}
\begin{cases}
\FDyn_1(\xx_1, \yy_1, \xx_2, \yy_2)
&:=\quad \xx_1 + \Pos{\yy_1} + \Pos{(\Pos{\yy_2} - \zz_1)},\\
\FDyn_2(\xx_1, \yy_1, \xx_2, \yy_2)
&:=\quad \yy_2 - \Pos\zz_1,\\
\FDyn_3(\xx_1, \yy_1, \xx_2, \yy_2)
&:=\quad  \xx_2 + \Neg{\yy_2} + \Neg{(\Neg{\yy_1} + \zz_1)},\\
\FDyn_4(\xx_1, \yy_1, \xx_2, \yy_2)
&:=\quad  \yy_1 + \Pos\zz_1,\\
\GDyn_1(\xx_1, \yy_1, \xx_2, \yy_2)
&:=\quad \xx_1 - \Pos{\yy_1} -  \Pos{(\Pos{\yy_2} + \zz_2)},\\
\GDyn_2(\xx_1, \yy_1, \xx_2, \yy_2)
&:=\quad  \yy_2 + \Neg\zz_2,\\
\GDyn_3(\xx_1, \yy_1, \xx_2, \yy_2)
&:=\quad  \xx_2 - \Neg{\yy_2} - \Neg{(\Neg{\yy_1} -\zz_2)},\\
\GDyn_4(\xx_1, \yy_1, \xx_2, \yy_2)
&:=\quad  \yy_1 - \Neg\zz_2,
\end{cases}
\end{equation}
where we put 
$\zz_1 := \xx_1 - \Neg{\yy_1} - \xx_2 + \Pos{\yy_2}$
and $\zz_2 := \xx_1 + \Neg{\yy_1} - \xx_2 - \Pos{\yy_2}$. 

$(ii)$ For $(\aa_1, \bb_1, ..., \aa_\nn, \bb_\nn)$
in~$\Int^{2n}$, put $(\aa_1, \bb_1, ..., \aa_\nn, \bb_\nn) \act
\s_\ii^\ee =  (\aa'_1, \bb'_1, ..., \aa'_\nn, \bb'_\nn)$ with
$\aa'_\kk = \aa_\kk$ and $\bb'_\kk=\bb_\kk$
for $\kk\not=\ii,\ii+1$, and
\begin{equation*}
(\aa'_\ii, \bb'_\ii, \aa'_{\ii+1}, \bb'_{\ii+1}) =
\begin{cases}
\FDyn(\aa_\ii, \bb_\ii, \aa_{\ii+1}, \bb_{\ii+1})
&\text{for $\ee=+1$},\\
\GDyn(\aa_\ii, \bb_\ii, \aa_{\ii+1}, \bb_{\ii+1})
&\text{for $\ee=-1$}.\\
\end{cases}
\end{equation*}
Then, for $\ww$ an $n$-braid word, we recursively define
\begin{equation*}
(\aa_1, \bb_1, ..., \aa_\nn, \bb_\nn) \act \ww =
\begin{cases}
(\aa_1, \bb_1, ..., \aa_\nn, \bb_\nn)
&\text{for $\ww= \e$},\\
((\aa_1, \bb_1, ..., \aa_\nn, \bb_\nn) \act \ww') \act \s_\ii^\ee
&\text{for $\ww = \ww'\s_\ii^\ee$}.
\end{cases}
\end{equation*}
The {\it Dynnikov coordinates} of~$\ww$
are defined to be the sequence $(0,1,0,1,..., 0,1) \act\ww$.
\end{defi}

\begin{theo}
\cite[Propositions~8.5.3 and~8.5.4]{Dgr}
\label{T:Dynnikov}
The Dynnikov coordinates of an $n$-braid word~$\ww$ characterize
the braid~$\cl\ww$ represented by~$\ww$: the coordinates
of~$\ww$ and~$\ww'$ are equal if and only if $\cl\ww=\cl{\ww'}$ holds.
\end{theo}

We deduce still another solution to the braid isotopy problem:

\begin{coro}
An $n$-strand braid word represents~$1$ in~$B_\nn$ if
and only if its Dynnikov coordinates are~$(0,1,0,1,..., 0,1)$.
\end{coro}

\begin{exam}
\label{X:Dynnikov}
Consider $\ww_0 = \mathtt{aBabacABABAbbCB}$ once more. By
applying the formulae of~\eqref{E:Dynnikov} inductively, we
compute the Dynnikov coordinates of the successive prefixes
of~$\ww_0$:

$0: \e \gives (0, 1, 0, 1, 0, 1, 0, 1)$

$1: \mathtt{a} \gives (1, 0, 0, 2, 0, 1, 0, 1)$

$2: \mathtt{aB} \gives (1, 0, -2, 0, 0, 3, 0, 1)$

$3: \mathtt{aBa} \gives (1, -3, -2, 3, 0, 3, 0, 1)$

$4: \mathtt{aBab} \gives (1, -3, 3, 2, 0, 4, 0, 1)$

$5: \mathtt{aBaba} \gives (1, -1, 3, 0, 0, 4, 0, 1)$

\ldots

$14: \mathtt{aBabacABABAbbC} \gives (1, -7, 5, -1, -7, 4, 0, 8)$

$15: \mathtt{aBabacABABAbbCB} \gives (1, -7, -6, 4, 1, -1, 0, 8)$.

\noindent The latter coordinates are not $(0,1,0,1,0,1,0,1)$, so we conclude that $\ww_0$ does not
represent~$1$ in~$B_4$.
\end{exam}

\subsubsection{Explanation}

At first, the formulae~\eqref{E:Dynnikov} seem quite mysterious. Actually,
there is no miracle here, but a very clever use of the simple formula
\begin{equation}
\label{E:Flip}
\xx+\xx'=\max(\xx_1+\xx_3,\xx_2 + \xx_4)
\end{equation}
that compares the number of intersections of a family of curves with
two triangulations obtained one from the other by switching one
diagonal in a quadrilateral (flip transformation). 

The framework consists in considering an $n$-braid as the isotopy
class of a homeomorphism of a disk with $n$~punctures: then the
generator~$\ss\ii$ corresponds to the homeomorphism that
exchanges the $\ii$th and $(\ii+1)$st punctures by a half-turn.
Dynnikov's idea is to let the braid act on a family of curves in the
punctured disk, and to count their intersections with a fixed
triangulation, or, equivalently, to let the braid act on
the triangulation and count its intersections with a fixed family of
curves. The action of~$\ss\ii$ on the chosen triangulation can then
be decomposed into the composition of four flips, and
applying~\eqref{E:Flip} repeatedly leads to the mysterious
formulae~\eqref{E:Dynnikov}. 

\subsubsection{Discussion}

The remarkable point about the Dynnikov coordinates is that they
involve the semiring $(\Int, \max, +, 0)$, which explains their
efficiency. Multiplying by one generator~$\ss\ii$ can only increase
the size of the coordinates by one unit, whereas similar formulae in
the ring $(\Int, +, \times, 1)$ would double
the size in the worst case. It follows that the solution to the braid word
problem given by Theorem~\ref{T:Dynnikov} has a linear space
complexity, and a quadratic time complexity.

It is not so easy to compare the solution based on the Dynnikov
coordinates with the other solutions to the braid isotopy problem,
because its practical efficiency much depends on the way large
integer arithmetic is implemented---large integers do appear
when long braid words are considered, typically length~$O(\ell)$ binary
integers for a length~$\ell$ braid word representing a pseudo-Anosov 
braid. However, the only
arithmetic operations involved are addition and maximum, and
both can be implemented very efficiently and easily. No statistical
study has been completed so far, and it would be desirable to
compare Dynnikov's method with handle reduction. The only weak
point of the former is that, at the moment, it is purely
incremental: the only available formulae correspond to multiplying by
one single letter~$\ss\ii$ or~$\sss\ii$, so, in particular, no
divide-and-conquer variant exists.

\end{document}